\documentclass[a4paper, 11pt, twoside]{article}
\usepackage[a4paper,margin=2.5cm]{geometry}
\usepackage[utf8]{inputenc}
\usepackage{lmodern}
\usepackage{amssymb}
\usepackage[english]{babel}
\usepackage{graphicx} 
\usepackage{epsfig}
\usepackage{epstopdf}
\usepackage{amsfonts} 
\usepackage{amsmath} 
\usepackage{braket}
\usepackage{amsthm} 
\usepackage{chngpage}
\usepackage[usenames,dvipsnames]{xcolor}
\usepackage{mathrsfs}
\usepackage{enumitem}
\usepackage{comment}
\usepackage{chngcntr}
\usepackage[final]{microtype}
\usepackage{todonotes}
\usepackage{array}
\usepackage{caption}
\usepackage{subcaption}
\usepackage[T1]{fontenc}
\usepackage{listings}
\usepackage{tikz}
\usepackage{tikz-cd}
\usepackage{hyperref}
\hypersetup{
    colorlinks,
    linkcolor={red!50!black},
    citecolor={blue!50!black},
    urlcolor={blue!80!black}
}
\usepackage{cleveref}
\lstdefinelanguage{Julia}
  {morekeywords={abstract,break,case,catch,const,continue,do,else,elseif,
      end,export,false,for,function,immutable,import,importall,if,in,
      macro,module,otherwise,quote,return,switch,true,try,type,typealias,
      using,while},
   sensitive=true,
   morecomment=[l]\#,
   morecomment=[n]{\#=}{=\#},
   morestring=[s]{"}{"},
   morestring=[m]{'}{'},
}[keywords,comments,strings]

\newcommand{\C}{\mathbb{C}}

\newcommand{\R}{\mathbb{R}}

\newcommand{\LA}{\mathcal{L}_A}
\newcommand{\LAb}{\mathcal{L}_{A,b}}

 \newcommand{\Row}{{\text{\rm Row}}}
\newcommand{\ML}{{\mathcal L}}
\newcommand{\Z}{\mathbb Z}

\newcommand{\rank}{\text{\rm rank}}

\newcommand{\azero}{\mathbf{0}}

\theoremstyle{plain}
	\newtheorem{theorem}{Theorem}[section] 	
	\newtheorem{lemma}[theorem]{Lemma}
	\newtheorem{corollary}[theorem]{Corollary}
	\newtheorem{proposition}[theorem]{Proposition}
    
     \newtheorem{definition}[theorem]{Definition}
     \newtheorem{example}[theorem]{Example}
	\newtheorem{remark}[theorem]{Remark}

\setlist{itemsep=2pt, topsep=3pt, parsep=0pt}

\let\OLDthebibliography\thebibliography
\renewcommand\thebibliography[1]{
  \OLDthebibliography{#1}
  \setlength{\parskip}{0.4pt}
  \setlength{\itemsep}{2.8pt plus 0.0ex}
}

\title{\bf Intersections, Equations, and Singularities\\ of Lissajous Varieties}
\author{ Francesco Maria Mascarin \thanks{MPI-MiS Leipzig, and CWI, Amsterdam}
}
\date{}
 
\begin{document}
 \maketitle

\begin{abstract}
\noindent Lissajous varieties are affine algebraic varieties parametrized by sine and cosine functions. They generalize algebraic Lissajous curves in the plane to higher dimensions by applying a component-wise trigonometric map to a rational affine subspace and are motivated by applications in dynamical systems (Kuramoto model) and semidefinite optimization (graph elliptopes). We establish several algebraic and geometric properties of Lissajous varieties. In particular, we show that intersections of Lissajous varieties are finite unions of such varieties and describe them explicitly. We investigate algebraic properties of the defining polynomials of Lissajous hypersurfaces and describe the singular locus of Lissajous varieties.
\end{abstract}

\section{Introduction}\label{sec:Introduction}
Let $A\in \Z^{d\times n}$ ($d<n$) be an integer matrix of rank $d$, $b\in \C^n$ a complex vector, and $L_{A,b}$ the affine complex subspace defined by $L_{A,b} = {\rm Row}(A) -  \frac{b\pi}{2} = \{ x -  \frac{b \pi}{2}\, : \, x \in {\rm Row}(A)\} \subseteq \mathbb{C}^n$.
The \emph{Lissajous variety} of the couple $(A,b)$ is 
\begin{equation}\label{eq:defLAb}
    \LAb = \cos(L_{A,b}),
\end{equation}
where $\cos$ indicates the complex cosine map $\cos(z)=\frac{1}{2}(e^{iz}+e^{-iz})$ applied component-wise, with $i=\sqrt{-1}$ the imaginary unit.
For convenience, we sometimes write $\LA = {\cal L}_{A,{\bf 0}}$ for the case $b={\bf 0}$.

Lissajous curves in the plane are parametric curves describing the superposition of two perpendicular oscillators whose angular frequencies ratio is rational \cite{greenslade1993all}.
They were named after Jules Antoine Lissajous and their applications span from engineering (including electronics, where they can be generated using an oscilloscope, and robotics \cite{borkar2020application}) to quantum mechanics \cite{jaczewski2024lissajous,russo2026lissajous}.

Lissajous varieties were first introduced in \cite{mascarin2025lissajousvarieties} and generalize algebraic Lissajous planar curves to higher dimensions.
They are motivated by applications in dynamical systems and mathematical optimization.
In dynamical systems, Lissajous varieties provide an algebraic encoding of the equilibria of the Kuramoto model, which describes the dynamics of a system of coupled phase oscillators \cite{kuramoto1975self}, \cite[Section 4]{mascarin2025lissajousvarieties}.
In mathematical optimization, Lissajous varieties (for $b={\bf 0}$) arise in the description of the algebraic boundary of graph elliptopes: convex semi-algebraic sets of particular relevance to semidefinite programming \cite{Laurent2026AlgebraicBoundary}.
Moreover, they represent a special instance of Chebyshev varieties, studied in \cite{BelAfiaMeroniTelen2025Chebyshev}.

\begin{example}\label{ex:ex1.1}
    The \emph{lemniscate} (eight- or infinity-figure) parametrized by $(\sin(t),\sin(2t))$ (see Figure \ref{fig:ex1.1}, on the left) is an example of Lissajous planar curve with $A=\left(\begin{smallmatrix}1&2\end{smallmatrix}\right)$ and $b=\bf 1$.

    Lissajous varieties with $b={\bf 1}$ arising from the incidence matrix $A(G)$ of a graph are of particular interest for encoding the equilibria of the Kuramoto model algebraically. 
    Figure \ref{fig:ex1.1} on the right displays the Lissajous variety ${\cal L}_{A,{\bf 1}}$ for $A=A(C_3)$, the incidence matrix of the triangle graph $C_3$.
    \begin{figure}
    \centering
    \includegraphics[height = 5cm]{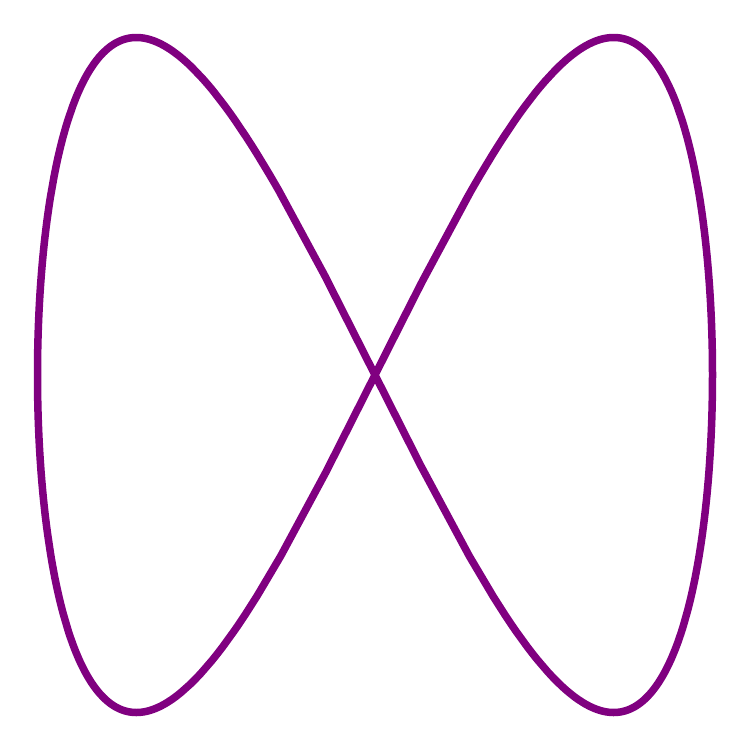} \quad  \quad \quad \quad 
    \includegraphics[height = 5cm]{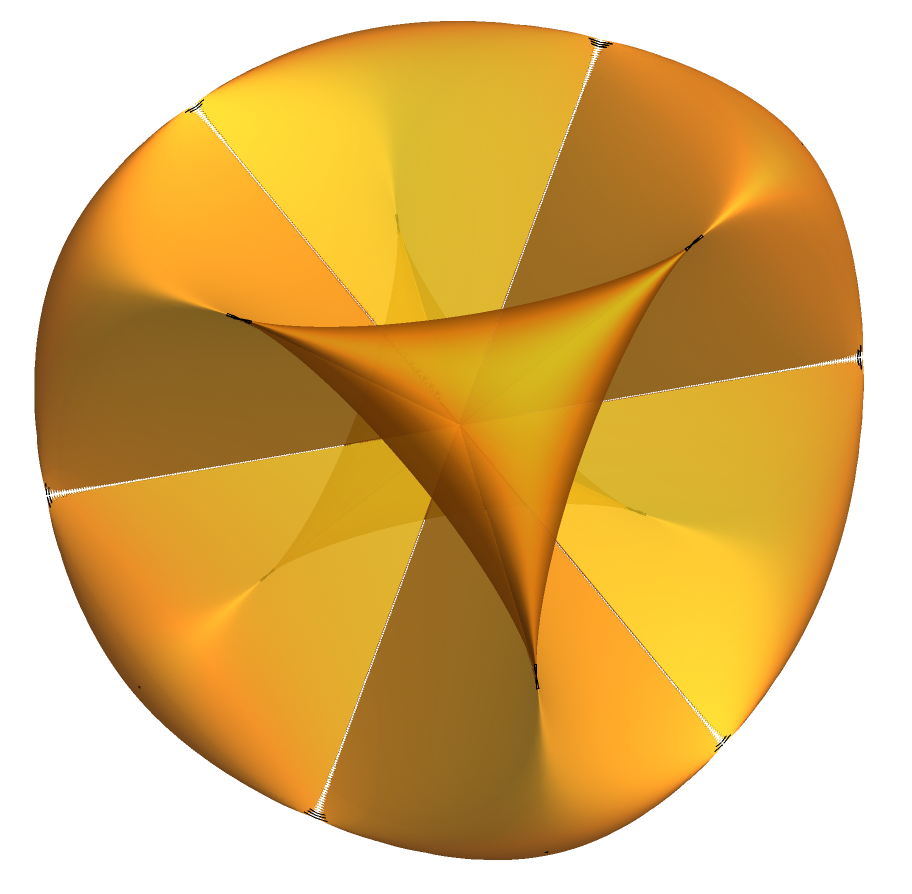}
    \caption{The lemniscate (left) and ${\cal L}_{A,{\bf 1}}$ (right) from Ex. \ref{ex:ex1.1} are examples of Lissajous varieties.}
    \label{fig:ex1.1}
    \end{figure}
\end{example}

Building on recent advances related to Lissajous varieties, this paper investigates intersections, defining equations of Lissajous hypersurfaces, and singularities. 

We show that Lissajous varieties are closed under taking intersections, in the sense that the intersection of any two Lissajous varieties is a finite union of Lissajous varieties, which we describe explicitly (Theorem \ref{thm:GeneralLissajousIntersection}).
Algebraic varieties of codimension one are called hypersurfaces and can be described by a single polynomial equation. 
We describe the defining polynomial of a Lissajous hypersurface using the determinantal representation proposed in \cite[Section 3]{mascarin2025lissajousvarieties}.
We show that this polynomial has constant leading coefficient when viewed as a univariate polynomial in any given variable.
For generic $b$, we provide a formula for its leading coefficient, degree and roots (Theorem \ref{thm:LissHypersurfaces}).
Finally, we study the singular locus of Lissajous varieties, establishing a sandwich-type double inclusion in terms of the component-wise cosine map applied to subsets of $L_{A,b}$.
This leads to the construction of two classes of subsets, $S_{A,b},W_{A,b}\subset L_{A,b}$, arising from rank-deficiency points of the parametrization and self-intersections, respectively.
In particular, we show that for curves ($d=1$, $n$ arbitrary) and hypersurfaces ($d=n-1$, $n\geq4$) the sets in $W_{A,b}$ alone suffice to characterize the singularities.
We illustrate these results in an example.

\begin{example}\label{ex:intro}
Let $A=\left ( \begin{smallmatrix} 3&-1&0\\5&0&-1\end{smallmatrix} \right)$, so that ${\cal L}_A=\{(\cos t_1,\cos t_2,\cos t_3)\;:\; t\cdot v=0\}\subset\mathbb{C}^3$, with $v=(1,3,5)$.
Its defining polynomial is
\begin{align*}
    \Gamma_A=&x_1^2 + 9 x_2^2 - 24 x_2^4 + 16 x_2^6 + 30 x_1 x_2 x_3 - 
 40 x_1 x_2^3 x_3 + 25 x_3^2 - 120 x_1 x_2 x_3^3 + 160 x_1 x_2^3 x_3^3\\
 &- 200 x_3^4 + 96 x_1 x_2 x_3^5 - 128 x_1 x_2^3 x_3^5 + 560 x_3^6 - 640 x_3^8 + 256 x_3^{10} -1.
\end{align*}
We can write ${\cal L}_A$ as the projection of the intersection of two Lissajous varieties in $\mathbb{C}^4$. 
Let
\[
A_1=\left ( \begin{matrix} 1&0&0&1\\0&1&0&3\\0&0&1&0\end{matrix} \right) \quad \text{and} \quad A_2=\left ( \begin{matrix} 1&0&0&0\\0&1&0&0\\0&0&1&-5\end{matrix} \right).
\]
Their row spaces intersect in $\{t_4=t_1+3t_2=-5t_3\}$ and ${\cal L}_{A_1}\cap{\cal L}_{A_2}$ is a single Lissajous variety (see Corollary~\ref{cor:intersection}).
Therefore ${\cal L}_A=\pi_{x_4}({\cal L}_{A_1}\cap{\cal L}_{A_2})$, where $\pi_{x_4}\colon\C^4\to\C^3$ is the projection onto the variables $(x_1,x_2,x_3)$.
Let $T_k$ be the $k$-th Chebyshev polynomial of the first kind, so that
$T_3(x)=4x^3-3x$ and $T_5(x)=16x^5-20x^3+5x$. 
The defining polynomials of ${\cal L}_{A_1}$ and ${\cal L}_{A_2}$ are
\[
\Gamma_{A_1}=x_1^2+x_4^2+T_3(x_2)^2-2x_1x_4T_3(x_2)-1,\qquad
\Gamma_{A_2}=x_4-T_5(x_3),
\]
respectively.
The polynomial $\Gamma_A$ is obtained as their resultant with respect to $x_4$:
\[
\Gamma_A=\operatorname{Res}_{x_4}(\Gamma_{A_1},\Gamma_{A_2})
=\Gamma_{A_1}\big|_{x_4=T_5(x_3)}
=x_1^2+T_3(x_2)^2+T_5(x_3)^2-2\,x_1T_3(x_2)T_5(x_3)-1 .
\]
The polynomial $\Gamma_{A_1}$ can itself be obtained as a resultant, of the Cayley cubic $x_1^2+x_4^2+x_5^2-2x_1x_4x_5-1$ and $x_5-T_3(x_2)$, eliminating the additional variable $x_5$.
Resultant representations of Lissajous hypersurfaces are discussed in full detail in Section \ref{sec:hypersurfaces}.

Theorem~\ref{thm:SingSandwich} gives a list of candidate components for ${\rm Sing}(\mathcal L_A)$, shown in Figure~\ref{fig:intro}.
The purple curves come from self-intersections of the parametrization.
They are the images under $\cos$ of the lines $\{t\in \Row(A):\ t_2\in\tfrac{\pi}{3}\mathbb Z\setminus\pi\mathbb Z\}$ and $\{t\in \Row(A):\ t_3\in\tfrac{\pi}{5}\mathbb Z\setminus\pi\mathbb Z\}$.
These lines satisfy properties (P1) and (P2) of Theorem \ref{thm:SingSandwich}, so by the first inclusion of the theorem their images lie in ${\rm Sing}(\mathcal L_A)$.
The blue points come from rank-deficiency points of the parametrization, i.e., points where its Jacobian drops rank:
$(1,1,1),\,(1,-1,-1),\,(-1,1,-1),\,(-1,-1,1)$.
They are singular because $\nabla\Gamma_A$ vanishes there.
By the second inclusion of Theorem~\ref{thm:SingSandwich}, there are no further singularities, since the remaining rank-deficiency points lie on the purple curves.
Hence ${\rm Sing}(\mathcal L_A)$ is the union of the six curves and the four points.
\begin{figure}
    \centering
    \includegraphics[height = 8cm]{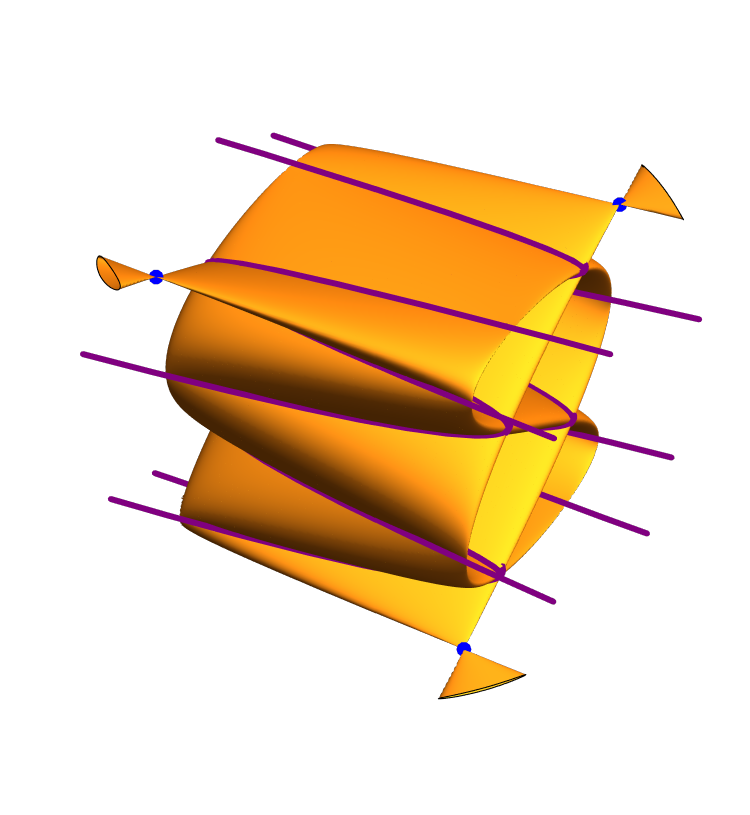} 
    \includegraphics[height = 8cm]{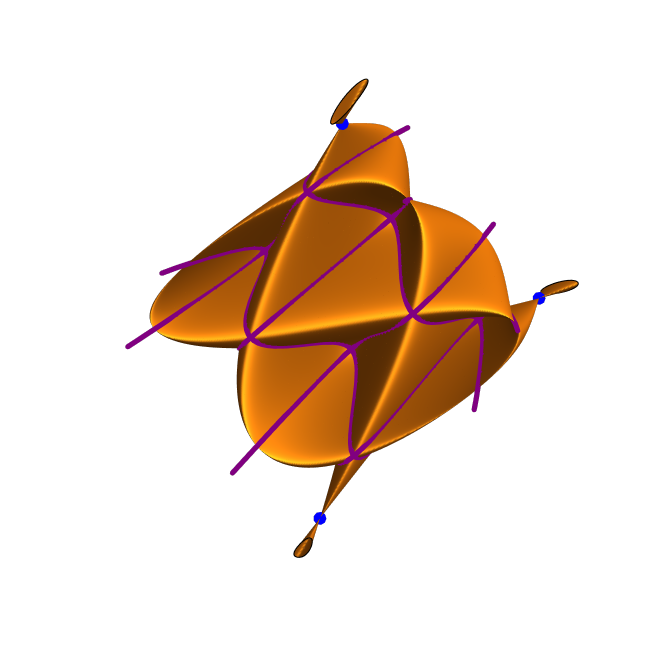}
    \caption{The Lissajous variety $\LA$ (gold) and its singular locus (purple, blue) from Example \ref{ex:intro}.}
    \label{fig:intro}
\end{figure}
\end{example}

\paragraph{Related work.} Lissajous varieties were first introduced in \cite{mascarin2025lissajousvarieties}, where their most relevant algebro-geometric properties were established.
This work built on an earlier study of algebraic varieties parametrized by Chebyshev polynomials \cite{BelAfiaMeroniTelen2025Chebyshev}; allowing a sine-based parametrization as well makes Lissajous varieties suitable for studying the equilibria of the Kuramoto model.
An alternative algebraic encoding of these equilibria was introduced in \cite{harrington2023kuramoto}, where the authors study Kuramoto varieties.

\emph{Graph elliptopes} correspond to sets of partially specified correlation matrices that admit a positive semidefinite (PSD) completion.
They are convex semi-algebraic sets and projected spectrahedra.
In \cite{Laurent2026AlgebraicBoundary}, it was shown that the non-trivial components of the algebraic boundary of graph elliptopes are Lissajous hypersurfaces constructed from cycles in the graph.
Their defining polynomials, called \emph{cycle polynomials}, were introduced in \cite{Sturmfels2010MultivariateGaussians}.
The same type of polynomials recently arose independently as components of \emph{matrix completion discriminants} for cycle completable graphs in \cite{Sekiguchi2026completion}. 
There they are called the \emph{iterated resultant} ${\rm ItrRes}(C_n)$, underlying a property of Lissajous hypersurfaces that will be exploited in Section \ref{sec:hypersurfaces}.

\paragraph{Organization of the paper.} In Section \ref{sec:intersections}, we show that the intersection of two Lissajous varieties is a finite union of Lissajous varieties (Theorem \ref{thm:GeneralLissajousIntersection}).
For $b= {\bf 0}$, we also show an instance in which this union collapses to a single variety.

In Section \ref{sec:hypersurfaces}, we describe the defining equation of Lissajous hypersurfaces. 
The main result here is Theorem \ref{thm:LissHypersurfaces}, which describes the polynomial for generic $b$ as a univariate polynomial in each of its variables.
For the special case $b={\bf 0}$, we provide a procedure for obtaining a resultant representation. 

In Section \ref{sec:singularities}, we describe the singular locus of Lissajous varieties, establishing a sandwich-type inclusion (Theorem \ref{thm:SingSandwich}).
In particular, we completely characterize the singularities of Lissajous curves in any dimension and of Lissajous hypersurfaces when $n>4$.

\section{Notation and Preliminaries}

Choosing coordinates on ${\rm Row}(A)$, the Lissajous variety ${\cal L}_{A,b}$ can be parametrized via the map: 
\begin{equation}\label{eq:parametrization}
    \phi_{A,b} \, : \, \mathbb{C}^d \longrightarrow \mathbb{C}^n, \quad t = (t_1, \ldots, t_d) \longmapsto ( {\rm cos}(a_1 \cdot t - b_1 \tfrac{\pi}{2}), \ldots, {\rm cos}(a_n \cdot t - b_n \tfrac{\pi}{2} )). 
\end{equation}
Here $a_j \in \mathbb{Z}^d$ is the $j$-th column of $A$ and $a_j \cdot t$ is the standard dot product. 
We clearly have ${\cal L}_{A,b} = {\rm im} \, \phi_{A,b}$. 
In \cite{mascarin2025lissajousvarieties}, it was shown that the Lissajous variety of $(A,b)$ is an irreducible affine variety of dimension $\rank(A)$, whose degree is derived from the volume of a polytope.
Let $A\in\Z^{d\times n}$ have rank $d$.
Let $\Z A\subseteq \Z^d$ be the sublattice of $\mathbb{Z}^d$ generated by $a_1,\dotsc,a_n\in\Z^d$, and assume that $\Z A = \Z^d$. 
This is not restrictive because after a linear change of coordinates on $\mathbb{C}^d$ one can replace $A$ by a matrix $A'$ with $\mathbb{Z}A'=\mathbb{Z}^d$ and $\Row(A)=\Row(A')$, which gives the same Lissajous variety.
\begin{theorem}[{\cite[Theorem 1.4]{mascarin2025lissajousvarieties}}]\label{thm:mainintro}
Let $A \in \mathbb{Z}^{d \times n}$ be such that ${\rm rank}(A) = d$ and $\mathbb{Z}A = \mathbb{Z}^d$. 
Let $P_A \subset \mathbb{R}^d$ be the polytope obtained as the convex hull of the lattice points $\{ \pm a_1, \ldots, \pm a_n \} \subset \mathbb{Z}^d$. 
For any $b \in \mathbb{C}^n$, 
\[
{\cal L}_{A,b} = \cos(L_{A,b}) = {\rm im}(\phi_{A,b})
\]
is an irreducible variety of dimension $d$. 
Moreover, for generic $b \in \mathbb{C}^n$, its degree is $\deg ({\cal L}_{A,b}) = 2^{-{\rm CL}_A} \, d! \, {\rm vol}(P_A)$, where ${\rm vol}(\cdot)$ denotes the euclidean volume and ${\rm CL}_A$ is the number of zero entries of a generic vector in the kernel of $A: \mathbb{R}^n \rightarrow \mathbb{R}^d$. 
\end{theorem}

For $z\in(\mathbb{C}^*)^k$ and $m\in\mathbb{Z}^k$ we write $z^{m}=z_1^{m_1}\cdots z_k^{m_k}$ for the multi-index power (for $z\in\mathbb{C}^k$, we use $m\in\mathbb{N}^k$ instead). 
For $z,z'\in\mathbb{C}^n$ we write $z\circ z'=(z_1z'_1,\dotsc,z_nz'_n)$ for the Hadamard (entrywise) product.

The integer matrix $A$ defines an affine toric variety $Y_A \subseteq \mathbb{C}^n$.
This is obtained as the closure of the image of the Laurent monomial map $v \mapsto (v^{a_1}, \ldots, v^{a_n})$, $v \in (\mathbb{C}^*)^d$ \cite{CoxLittleSchenck2011,telen_2026_22806644}. 
Its ideal $I(Y_A)$ is generated by binomials of the form $y^{u} - y^{w}\in \mathbb{C}[y_1, \ldots, y_n]$, with $u, w \in \mathbb{N}^n$ and $A (u-w) = 0$ (\cite[Theorem~1.3.2]{telen_2026_22806644}).
A \emph{scaled} toric variety $Y_{A,\beta}$ is obtained from $A \in \mathbb{Z}^{d \times n}$, $\beta \in (\mathbb{C}^*)^n$ as follows:
\[ 
Y_{A,\beta} \, = \, \{ \beta\circ y=(\beta_1 \, y_1, \ldots, \beta_n \, y_n) \, : \,  y \in Y_A \} \, \subseteq \, \mathbb{C}^n,
\]
This variety is parametrized by $v \mapsto (\beta_1 v^{a_1}, \ldots, \beta_n v^{a_n})$. 
One checks that if $I(Y_{A})$ is generated by the binomials $y^{u_k} - y^{w_k}, k = 1, \ldots, r$, then $I(Y_{A,\beta}) = \langle \beta^{w_1}y^{u_1} - \beta^{u_1}y^{w_1}, \ldots, \beta^{w_r}y^{u_r}-\beta^{u_r}y^{w_r} \rangle$. 
In this paper, we will be mostly interested in toric varieties in the torus $Y_A\cap (\C^*)^n$ and their scaled counterpart $Y_{A,\beta}\cap (\C^*)^n$.
For this reason, we denote with $I_A$ (resp., $I_{A,\beta}$) the ideal of $Y_A\cap (\C^*)^n$ (resp., $Y_{A,\beta}\cap (\C^*)^n$).
To compute the generators of $I_{A,\beta}$ one computes a $\mathbb{Z}$-kernel of $A$. 
If this kernel is generated by $m_1, \ldots, m_r \in \mathbb{Z}^n$, then the ideal is generated by the binomials $y^{m_k} - \beta^{m_k}, k = 1, \ldots, r$ (and so for $I_A$, with $\beta=\bf 1$) \cite[Prop.~1.2.22, Thm.~1.3.2]{telen_2026_22806644}.
The motivation for this framework will become clear in Section \ref{sec:intersections} and \ref{sec:hypersurfaces}, and relies on the following result from \cite{mascarin2025lissajousvarieties}, which shows that Lissajous varieties can be obtained as projections of incident varieties that are isomorphic to scaled toric varieties in the torus.

\begin{theorem}[{\cite[Theorem~2.2]{mascarin2025lissajousvarieties}}]\label{thm:projY_Ab}
    The Lissajous variety ${\cal L}_{A,b} \subseteq \mathbb{C}^n$ is the image of the variety
    \[ 
    {\cal Y}_{A,b} \, = \, \{(x,y) \in \mathbb{C}_x^{n}\times \mathbb{C}_y^{n} \, : \, y \in Y_{A,\beta}\cap(\mathbb{C}^*)^n \text{ and } x_j \, = \, \tfrac{1}{2}(y_j + y_j^{-1}) \text{ for } j = 1, \ldots,n \}, \]
    under the coordinate projection $\pi_{A,b}: {\cal Y}_{A,b}  \rightarrow \mathbb{C}^n$. 
    Here $\beta = (e^{-i b_1 \tfrac{\pi}{2}}, \ldots, e^{-i b_n \tfrac{\pi}{2}}) \in (\mathbb{C}^*)^n$.
\end{theorem}

In particular, from this result one can derive a formula for the degree of a Lissajous variety $\LAb$ that holds for any $b\in\C^n$ (Theorem \cite[Theorem 2.3]{mascarin2025lissajousvarieties}):
\begin{equation}\label{eq:LissDeg}
    \deg ({\cal L}_{A,b}) \, = \,  \frac{ d! \, {\rm vol}(P_A)}{\deg \pi_{A,b} \cdot [\mathbb{Z}^d : \mathbb{Z} A]}
\end{equation}

\section{Intersections of Lissajous Varieties}\label{sec:intersections}
In this section we show that Lissajous varieties are closed under taking  intersections. That is, the intersection of Lissajous varieties is a union of Lissajous varieties.

We proceed in three steps. 
First, we prove the statement for $\LA\cap{\cal L}_B$ under the assumption that the sum of the $\Z$-kernels of $A$ and $B$, after allowing sign changes in some columns of $B$, generates the full lattice in its real span.
Second, we remove the extra hypothesis on the lattice index.
Finally, we extend the result to the general case of ${\cal L}_{A_1,b_1}\cap{\cal L}_{A_2,b_2}$, allowing non-zero $b_1$ and $b_2$.

Let $A\in\mathbb{Z}^{d_1\times n}$ and $B\in\mathbb{Z}^{d_2\times n}$. 
Set $C\in\mathbb{Z}^{d_3\times n}$ such that ${\rm Row}(C)={\rm Row}(A)\cap{\rm Row}(B)$.
Since a Lissajous variety only depends on the row span of a matrix, we write
${\cal L}_{A\cap B,b}={\cal L}_{C,b}$.

We first consider the case $b=\mathbf 0$.
Using the trigonometric parametrization, one observes that the inclusion ${\cal L}_{A\cap B}\subseteq \LA \cap {\cal L}_B$ always holds.
Indeed, if $x\in \cos(\Row (C))$, with ${\rm Row}(C)={\rm Row}(A)\cap{\rm Row}(B)$, then clearly $x\in \cos(\Row (A))\cap \cos(\Row (B))$.
In general, the inclusion can be strict.
For instance, consider $A=\left(\begin{smallmatrix}1&0&-1\\-1&1&0\end{smallmatrix}\right)$ and $B=\left(\begin{smallmatrix}1&2&0\\0&0&1\end{smallmatrix}\right)$ from Example \ref{ex:intersectionlatticeindex1}.
Their row spaces intersect in $\Row(A)\cap\Row(B)=\langle(1,2,-3)\rangle$, hence ${\cal L}_{A\cap B}=\{(\cos t,\cos 2t,\cos 3t)\,:\,t\in\C\}$.
However, also $\{(\cos t,\cos 2t,\cos t)\,:\,t\in\C\}\subseteq \LA \cap {\cal L}_B$.
Indeed, $(\cos t,\cos 2t,\cos t)=\cos(t,2t,t)=\cos(t,-2t,t)$, where $(t,2t,t)\in{\rm Row}(B)=\{z\in\C^3:z_2=2z_1\}$ and $(t,-2t,t)\in{\rm Row}(A)=\{z\in\C^3:z_1+z_2+z_3=0\}$.
To capture the extra components in the intersection, we need to introduce some notation.

Let $J\subset[n]$ be a set of indices. 
We write $A^{-J}$ for the matrix $A$ where all the columns $a_j$ indexed with some $j\in J$ are replaced by $-a_j$.
When $v\in \C^n$ is a vector, we write $v_{-J}$ for $v$ with the signs of the entries in $J$ flipped, and $v^{-J}$ for $v$ with entries $j\in J$ replaced by $v_j^{-1}$.

\begin{theorem}\label{thm:LissajousIntersection}
    Let $A\in\mathbb{Z}^{d_1\times n}$ and $B\in\mathbb{Z}^{d_2\times n}$. 
    Set $M^{(J)}={\rm ker}_{\mathbb{Z}}(A)+{\rm ker}_{\mathbb{Z}}(B^{-J})$, viewed as a sublattice of $\mathbb{Z}^n\cap ({\rm ker}_{\R}(A)+{\rm ker}_{\R}(B^{-J}))$.
    If the index $[\mathbb{Z}^n\cap\R M^{(J)}:\mathbb{Z}M^{(J)}]=1$ for all $J\subset[n]$, then
    \[
    {\cal L}_A\cap{\cal L}_B=\bigcup_{J\subset[n]}{\cal L}_{A\cap B^{-J}}
    \]
\end{theorem}

\begin{example}\label{ex:intersectionlatticeindex1}
Let $A=\left(\begin{smallmatrix}1&0&-1\\-1&1&0\end{smallmatrix}\right)$ and $B=\left(\begin{smallmatrix}1&2&0\\0&0&1\end{smallmatrix}\right)$, so that $\ker_{\Z}(A)=\langle u\rangle$ and $\ker_{\Z}(B)=\langle v\rangle$ with $u=(1,1,1)$ and $v=(2,-1,0)$.
Then ${\cal L}_A$ is the Cayley cubic surface with equation $x_1^2+x_2^2+x_3^2-2x_1x_2x_3-1=0$ and ${\cal L}_B$ is the parabolic cylinder given by $x_2=2x_1^2-1$.
For every $J\subset[3]$, $v_{-J}=\pm(2,-1,0)$ or $v_{-J}=\pm(2,1,0)$, where $v_{-J}$ is $v$ with the signs of the entries in $J$ changed.
The $2\times2$ submatrix of $(u\,|\,v_{-J})$ on the last two rows has determinant $\pm1$ in both cases, so it is invertible over $\Z$.
Hence, if $\lambda u+\mu v_{-J}\in\Z^3$ with $\lambda,\mu\in\R$, then $(\lambda,\mu)$ is obtained from the last two coordinates of $\lambda u+\mu v_{-J}$ by applying the integer inverse, so $\lambda,\mu\in\Z$. 
Thus the lattice index is $1$ for all $J\subset[3]$.
In the first case ${\rm Row}(A)\cap{\rm Row}(B^{-J})=\langle(1,2,-3)\rangle$, and in the second $\langle(1,-2,1)\rangle$. 
By \Cref{thm:LissajousIntersection},
\[
{\cal L}_A\cap{\cal L}_B=\{(\cos t,\cos 2t,\cos 3t)\,:\,t\in\C\}\ \cup\ \{(\cos t,\cos 2t,\cos t)\,:\,t\in\C\},
\]
where only the first curve is ${\cal L}_{A\cap B}$ (see Figure~\ref{fig:intersections3.2and3.7}, on the left).
\end{example}

Before providing the proof, we need a preliminary lemma.

\begin{lemma}\label{lemma:Yintersection}
    Let $A\in\mathbb{Z}^{d_1\times n}$, $B\in\mathbb{Z}^{d_2\times n}$ and $M={\rm ker}_{\mathbb{Z}}(A)+{\rm ker}_{\mathbb{Z}}(B)$.
    Let $C\in\mathbb{Z}^{d_3\times n}$ with
    ${\rm Row}(C)={\rm Row}(A)\cap{\rm Row}(B)$. 
    If the lattice index $[\mathbb{Z}^n\cap\R M:\mathbb{Z}M]=1$, then
    \[
    {\cal Y}_{C}={\cal Y}_{A}\cap {\cal Y}_{B},
    \]
    where ${\cal Y}_{A}={\cal Y}_{A,\bf 0}$, and ${\cal Y}_{A,\bf 0}$ is the incidence variety from Theorem \ref{thm:projY_Ab}.
\end{lemma}
\begin{proof}
    We first claim that if $[\mathbb{Z}^n\cap\R M:\mathbb{Z}M]=1$, then a $\mathbb{Z}$-basis for $M$ is a $\mathbb{Z}$-basis for ${\rm ker}_{\mathbb{Z}}(C)$. 
    Indeed, by definition, $[\mathbb{Z}^n\cap\R M:\mathbb{Z}M]=1$ if and only if a $\mathbb{Z}$-basis for $M$ is a $\mathbb{Z}$-basis for $\mathbb{Z}^n\cap\R M=\mathbb{Z}^n\cap ({\rm ker}_{\R}(A)+{\rm ker}_{\R}(B))$. 
    On the other hand, ${\rm ker}_{\R}(C)={\rm Row}(C)^{\perp}=({\rm Row}(A)\cap{\rm Row}(B))^{\perp}={\rm Row}(A)^\perp+{\rm Row}(B)^\perp={\rm ker}_{\R}(A)+{\rm ker}_{\R}(B)$.
    Hence, $\mathbb{Z}^n\cap\R M={\rm ker}_{\mathbb{Z}}(C)$.

    Let $I_A$ be the ideal of $Y_A\cap (\C^*)^n$. 
    Then, $I_A=\langle y^u-1 \,:\, u\in\ker_{\mathbb{Z}}(A)\rangle = \langle y^{u_j}-1\, :\, \{u_j\}_j\, \text{ is a $\mathbb{Z}$-basis for} \ker_{\mathbb{Z}}(A)\rangle$. 
    Set $I_{A+B}=\langle y^w-1\, :\, w \in \ker_{\mathbb{Z}}(A)+\ker_{\mathbb{Z}}(B)\rangle$. 
    It is easy to check that $I_{A+B}=I_A+I_B$.
    In particular, 
    \begin{align*}
        I_A+I_B=I_{A+B}=\langle y^{w_j}-1\, :\, \{w_j\}_j \text{ is a $\mathbb{Z}$-basis for } M = \ker_{\mathbb{Z}}(A)+\ker_{\mathbb{Z}}(B)\rangle &=\\
        =\langle y^{w_j}-1\, :\, \{w_j\}_j \text{ is a $\mathbb{Z}$-basis for} \ker_{\mathbb{Z}}(C)\rangle &= I_C.
    \end{align*}
    Passing to their vanishing set, this means that $(Y_A\cap Y_B)\cap (\C^*)^n=Y_C\cap (\C^*)^n$. 
    In particular, using the definition, ${\cal Y}_{A}\cap {\cal Y}_{B}={\cal Y}_{C}$.
\end{proof}
\begin{proof}[Proof of \Cref{thm:LissajousIntersection}]
    Let $\pi_x \colon \C_x^n\times \C_y^n\rightarrow\C_x^n$ be the projection onto the $x$ coordinates. 
    We first show that $\pi_x({\cal Y}_A)\cap\pi_x({\cal Y}_B)=\bigcup_{J\subset[n]}\pi_x({\cal Y}_{A}\cap {\cal Y}_{B^{-J}})$.

    We first show the inclusion ``$\subseteq$''.
    Let $x^*\in\pi_x({\cal Y}_A)\cap\pi_x({\cal Y}_B)$. Then there exist $y$ and $z$ in $\C^n$ such that $(x^*,y)\in {\cal Y}_A$ and $(x^*,z)\in {\cal Y}_B$. Write $J_{x^*}=\langle 2 x^*_1 - y_1 - y_1^{-1}, \ldots, 2 x^*_n - y_n - y_n^{-1}\rangle \subset \C[y^{\pm}]$. Its vanishing set consists of at most $2^n$ points of the form $V(J_{x^*})=\{(y_1^{\pm1},\dotsc,y_n^{\pm 1})\}$. 
    Then, by definition of ${\cal Y}_A$, $y\in Y_A\cap(\C^*)^n$ and $y\in V(J_{x^*})$. 
    Similarly, $z\in Y_B\cap(\C^*)^n$ and $z\in V(J_{x^*})$. 
    Since $V(J_{x^*})=\{(y_1^{\pm 1},\dotsc,y_n^{\pm 1})\}$, then, for all $j\in [n]$, $z_j$ is either equal to $y_j$ or to $y_j^{-1}$. 
    Let $J\subset [n]$ be such that $z_j=y_j^{-1}$ for all $j\in J$, and write $z^{-J}$ for $z$ where $z_j$ is replaced by $z_j^{-1}$ for all $j\in J$. 
    Clearly, $y = z^{-J}$ and $z^{-J}\in Y_{B^{-J}}$. 
    In this way, we have shown that there exists $y\in (Y_A\cap Y_{B^{-J}})\cap(\C^*)^n$ for some $J\subset [n]$ such that $y\in V(J_{x^*})$. 
    Then, $x^*\in\pi_x({\cal Y}_{A}\cap {\cal Y}_{B^{-J}})$ for some $J\subset [n]$, which proves the first inclusion.

    For the opposite inclusion, let $x^*\in\pi_x({\cal Y}_{A}\cap {\cal Y}_{B^{-J}})$ for some $J\subset [n]$. 
    As before, there exists $y\in (Y_A\cap Y_{B^{-J}})\cap(\C^*)^n$, for some $J\subset [n]$, such that $y\in V(J_{x^*})$. 
    Since $y\in Y_{B^{-J}}$, we have $y^{-J}=z\in Y_B$.  
    Hence, we have found $y$ and $z$ in $\C^n$ such that $(x^*,y)\in {\cal Y}_A$ and $(x^*,z)\in {\cal Y}_B$, proving the claim.

    From Lemma \ref{lemma:Yintersection} applied to the pair $(A,B^{-J})$, ${\cal Y}_{A}\cap {\cal Y}_{B^{-J}}={\cal Y}_{C^{(J)}}$ for some $C^{(J)}\in\mathbb{Z}^{d_3\times n}$ such that ${\rm Row}(C^{(J)})={\rm Row}(A)\cap{\rm Row}(B^{-J})$. This implies that
    \[
    {\cal L}_A\cap{\cal L}_B = \pi_x({\cal Y}_A)\cap\pi_x({\cal Y}_B)=\bigcup_{J\subset[n]}\pi_x({\cal Y}_{A}\cap {\cal Y}_{B^{-J}}) = \bigcup_{J\subset[n]}\pi_x({\cal Y}_{C^{(J)}}) = \bigcup_{J\subset[n]}{\cal L}_{A\cap B^{-J}}.
    \]
\end{proof}

In Example \ref{ex:intersectionlatticeindex1}, the defining equations of ${\rm Row}(A)$ and ${\rm Row}(B)$ share two variables, and the sign changes produce a second component. 
When the defining equations share at most one variable, the sign changes have no effect and the union in \Cref{thm:LissajousIntersection} collapses to a single Lissajous variety.

\begin{corollary}\label{cor:intersection}
    Let $A,B\in\mathbb{Z}^{n-1\times n}$ be full rank matrices such that the linear equations defining ${\rm Row}(A)$ and ${\rm Row}(B)$ use at most one common variable.    
    If the index $[\mathbb{Z}^n\cap\R M:\mathbb{Z}M]=1$, where $M={\rm ker}_{\mathbb{Z}}(A)+{\rm ker}_{\mathbb{Z}}(B)$, then
    \[
    {\cal L}_A \cap {\cal L}_B = {\cal L}_{A\cap B}.
    \]
\end{corollary}
\begin{proof}
    Let $u,v\in \mathbb{Z}^n$ be such that $\ker(A)=\langle u\rangle$ and $\ker(B)=\langle v\rangle$.
    If $u$ and $v$ are linearly dependent, then the claim is trivial, so we assume they are independent.
    The lattice $\ker_{\mathbb Z}(A)+\ker_{\mathbb Z}(B^{-J})=\langle u,v_{-J}\rangle$, where $v_{-J}$ is $v$ with the signs of the entries in $J$ changed. 
    Its index in $\mathbb Z^n\cap\R\langle u,v_{-J}\rangle$ is the greatest common divisor of the $2\times2$ minors of $(u\,|\,v_{-J})$ (one can see this using Smith normal forms, see below). 
    Since $u$ and $v$ share at most one index, each minor $u_lv_k-u_kv_l$ has at most one non-zero term, so the minors of $(u\,|\,v_{-J})$ agree with those of $(u\,|\,v)$ up to sign.
    Hence the gcd does not depend on $J$, and the index is $1$ for all $J$ by assumption.
    Then, by Theorem \ref{thm:LissajousIntersection}, we get ${\cal L}_A\cap{\cal L}_B = \bigcup_{J\subset[n]}{\cal L}_{A\cap B^{-J}}$.
    
    To prove the statement, it suffices to show that ${\cal L}_{A\cap B}={\cal L}_{A\cap B^{-J}}$ for all $J\subset[n]$.

    Up to relabeling the variables, we may assume
    \[
    u=(u_1,\dotsc,u_p,0,\dotsc,0),\qquad v=(0,\dotsc,0,v_p,\dotsc,v_n).
    \]
    If the supports are disjoint, we choose $p$ so that this still holds, with $u_p=0$ or $v_p=0$.
    
    Let $J\subset[n]$ and let us first suppose that $p\in J$. Let $C$ be such that ${\rm Row}(C)={\rm Row}(A)\cap{\rm Row}(B)$. Write $C'$ for the matrix obtained from $C$ by changing the sign of every column whose index $j$ lies in $\{p+1,\dotsc,n\}\setminus J$. 
    We claim that $\ker(C')=\ker(A)+\ker(B^{-J})$. 
    Let $w'\in\ker(A)+\ker(B^{-J})$. 
    Since $\ker(B^{-J})=\langle v_{-J}\rangle$, we can write $w'=\lambda u+\mu v_{-J}$ and
    \begin{align*}
    (C'w')_l&=\sum_{j=1}^{p-1}(C_{lj} \cdot \lambda u_j) + C_{lp}(\lambda u_p - \mu v_p) +\!\! \sum_{ \{p+1,\dotsc,n\}\setminus J}\! -C_{lj} \cdot \mu v_j +\!\! \sum_{J\cap\{p+1,\dotsc,n\}\! }C_{lj} \cdot(- \mu v_j)=\\
    &=\lambda \sum C_{lj}u_j - \mu \sum C_{lj}v_j = \lambda(Cu)_l-\mu(Cv)_l=0.
    \end{align*}
    Hence, $w'\in \ker(C')$. Since both $\ker(C')$ and $\ker(A)+\ker(B^{-J})$ have dimension $2$ as linear spaces, equality holds.

    Let now suppose $p\not\in J$. Let $C'$ be the matrix obtained from $C$ by changing the sign of all columns whose index $j$ lies in $\{1,\dotsc, p\}$ or $\{p+1,\dotsc,n\}\setminus J$. 
    Again, we claim that $\ker(C')=\ker(A)+\ker(B^{-J})$. 
    Let $w'= \lambda u+\mu v_{-J}\in\ker(A)+\ker(B^{-J})$. Then 
    \begin{align*}
    (C'w')_l&=\sum_{j=1}^{p-1}(-C_{lj} \cdot \lambda u_j) - C_{lp}(\lambda u_p + \mu v_p) +\!\! \sum_{ \{p+1,\dotsc,n\}\setminus J}\! -C_{lj} \cdot \mu v_j +\!\! \sum_{J\cap\{p+1,\dotsc,n\}}\! C_{lj} \cdot(- \mu v_j)=\\
    &=-\lambda \sum C_{lj}u_j - \mu \sum C_{lj}v_j = (-C(\lambda u+\mu v))_l=0.
    \end{align*}
    As before, we can conclude that $\ker(C')=\ker(A)+\ker(B^{-J})$.

    In both cases, by construction $\cos({\rm Row}(C))=\cos({\rm Row}(C'))$, that is ${\cal L}_{A\cap B}={\cal L}_{A\cap B^{-J}}$.
\end{proof}

Corollary \ref{cor:intersection} provides the underlying theory for the construction in Example \ref{ex:intro}. 
There, the kernels of $A_1$ and $A_2$ are generated by $(1,3,0,-1)$ and $(0,0,5,1)$, whose supports meet only in the index $4$, and the $2\times 2$ minors of the matrix with these two columns have greatest common divisor $1$. 
Thus, as in the proof of Corollary~\ref{cor:intersection}, the lattice index hypothesis is satisfied. 
Hence ${\cal L}_{A_1}\cap{\cal L}_{A_2}={\cal L}_{A_1\cap A_2}$, and projecting out the shared variable $x_4$ yields ${\cal L}_A$.

\begin{figure}
    \centering
    \includegraphics[height=5cm]{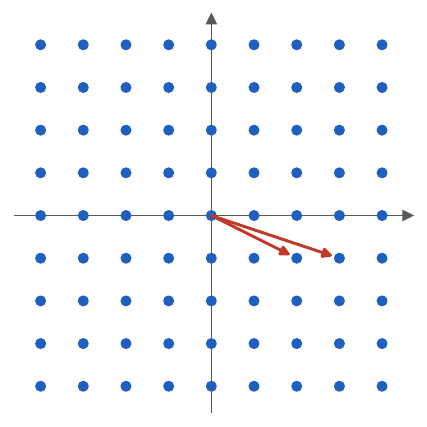}\qquad
    \includegraphics[height=5cm]{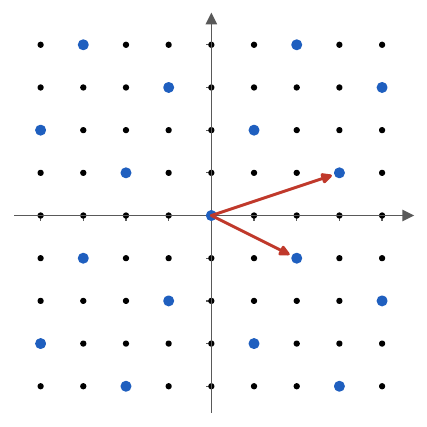}
    \caption{The lattice $M^{(J)}$ (blue) inside $\Z^2$ (black), with its generators (red), for $A=(1\ 2)$ and $B=(1\ 3)$: $J=\emptyset$, with index $1$ (left), and $J=\{2\}$, with index $5$ (right).}
    \label{fig:latticeindex}
\end{figure}

In general, the lattice index $[\mathbb{Z}^n\cap\R M^{(J)}:\mathbb{Z}M^{(J)}]$ may be greater than one for some $J\subset[n]$.
For instance, consider $A=\left(\begin{smallmatrix} 1&2 \end{smallmatrix}\right)$ and $B=\left(\begin{smallmatrix} 1&3 \end{smallmatrix}\right)$ from Example \ref{ex:intersection1213}.
The lattice $M$ is generated by $(2,-1)$ and $(3,-1)$ and has index $1$ in $\Z^2$.
On the other hand, for $J=\{2\}$, $M^{(J)}$ is generated by $(2,-1)$ and $(3,1)$.
Lattice points such as $(1,0)$ are not integer linear combinations of these two vectors, and indeed $M^{(J)}$ has index $5$ in $\Z^2$ (see Figure~\ref{fig:latticeindex}).
In this case \Cref{thm:LissajousIntersection} does not apply, and ${\cal L}_A\cap{\cal L}_B$ contains extra points coming from roots of unity.

To describe these extra components when the lattice index is greater than one, we use the Smith normal form of an integer matrix.
Let $M\in\mathbb{Z}^{n\times d}$ be an integer matrix. 
A \emph{Smith normal form} for $M$ is a diagonal matrix $S\in\mathbb{Z}^{n\times d}$ for which there exists $P\in\mathbb{Z}^{n\times n}$ and $Q\in\mathbb{Z}^{d\times d}$ with determinant $\pm 1$ (i.e., unimodular) such that $S=PMQ$. 
The diagonal entries $(s_1,\dotsc,s_d)$ of $S$ are called \emph{invariant factors}, and they are unique up to sign (reordering so that $s_i|s_{i+1}$).
The number of non-zero invariant factors reveals the rank of $M$, and, if $\rank(M)=d$, the lattice index of $\mathbb{Z}M$ in $\mathbb{Z}^n\cap\R M$ equals the absolute value of their product, namely $[\mathbb{Z}^n\cap\R M:\mathbb{Z}M]=|s_1\cdots s_d|$.

\begin{lemma}\label{lemma:toricintersection}
    Let $A\in\mathbb{Z}^{d_1\times n}$ and $B\in\mathbb{Z}^{d_2\times n}$, let $C$ be such that ${\rm Row}(C)={\rm Row}(A)\cap{\rm Row}(B)$, and let $d=\dim(\ker(A)+\ker(B))$.
    Let $M=(m_1|\dotsc|m_d)\in\mathbb{Z}^{n\times d}$, where $\{m_j\}_{j=1}^d$ is a $\mathbb{Z}$-basis of $\ker_{\mathbb{Z}}(A)+\ker_{\mathbb{Z}}(B)$, and let $S=PMQ$ be a Smith normal form of $M$ with invariant factors $s_1,\dotsc,s_d$.
    Set $Z=\{\zeta\in \C^n\, :\, \zeta_j^{s_j}=1 \text{ for } j\le d,\ \zeta_j=1 \text{ for } j>d\}$.
    Denote by $I_A$ the ideal of $Y_A\cap(\C^*)^n$ in $\C[y^{\pm1}]$, and by $I_{C,\xi}$ the ideal of $Y_{C,\xi}\cap(\C^*)^n$. Then
    \[
    I\big((Y_A\cap Y_B)\cap(\mathbb{C}^*)^n\big)=\bigcap_{\zeta \in Z} I_{C,\zeta^P}.
    \]
\end{lemma}

\begin{proof}
As in the proof of Lemma~\ref{lemma:Yintersection}, $I_A+I_B=\langle y^M-{\bf 1}\rangle$, where $y^M=(y^{m_1},\dotsc,y^{m_d})$.
Denote by $\varphi_M\colon (\C^*)^n\to(\C^*)^d$ the monomial map $y\mapsto y^M$.
Since $M=P^{-1}SQ^{-1}$, one checks that $\varphi_M=\varphi_{Q^{-1}}\circ\varphi_{S}\circ\varphi_{P^{-1}}$ (see, e.g., \cite[Ex. 1.2.11]{telen_2026_22806644}).
Consider the change of variables $z=\varphi_{P^{-1}}(y)$, i.e., $y=\varphi_P(z)=z^P$, and the $\mathbb{C}$-algebra isomorphism
\[
\varphi_P^*\colon \mathbb{C}[y^{\pm1}]\to\mathbb{C}[z^{\pm1}],\qquad f\mapsto f\circ\varphi_P,
\]
obtained by pulling back $\varphi_P$.
Then $\varphi_P^*(y^M)=z^{PM}=z^{SQ^{-1}}=(z^S)^{Q^{-1}}$, where $z^S=(z_1^{s_1},\dotsc,z_d^{s_d})$.
Since $\varphi_P^*$ is a ring isomorphism, it maps the ideal generated by $y^M-\mathbf 1$ to the ideal generated by $\varphi_P^*(y^M)-\mathbf 1$. 
Since $Q$ is unimodular, it holds that $\langle w^{Q^{-1}}-\mathbf 1\rangle=\langle w-\mathbf 1\rangle$. 
Therefore we get
\[
\varphi_P^*(I_A+I_B)
=\langle (z^S)^{Q^{-1}}-\mathbf 1\rangle
=\langle z^S-\mathbf 1\rangle\subset \C[z^{\pm1}],
\]
that is, $I_A+I_B=(\varphi_P^*)^{-1}\big(\langle z^S-\mathbf 1\rangle\big)$.

The ideal $\langle z^S-\mathbf 1\rangle\subset \C[z^{\pm 1}]$ is radical (one can see this through the Jacobian criterion), hence so is $I_A+I_B$, and therefore $I_A+I_B=I\big((Y_A\cap Y_B)\cap(\C^*)^n\big)$. 
Moreover,
\[
\langle z^S-\mathbf 1\rangle = \bigcap_{\zeta\in Z} \langle z_1-\zeta_1, \dotsc, z_d-\zeta_d \rangle.
\]

Since $MQ=P^{-1}S$ has columns $s_1(P^{-1})_1,\dotsc,s_d(P^{-1})_d$ and $Q$ is unimodular, these vectors form a $\mathbb Z$-basis of $\mathbb ZM$. 
As the columns of the unimodular matrix $P^{-1}$ form a $\mathbb Z$-basis of $\mathbb Z^n$, it follows that $(P^{-1})_1,\dotsc,(P^{-1})_d$ form a $\mathbb{Z}$-basis of $\mathbb{Z}^n\cap\R M$.
Since $\ker_{\mathbb Z}(C)=\mathbb Z^n\cap\R M$, these columns form a $\mathbb Z$-basis of $\ker_{\mathbb Z}(C)$, which is what identifies the ideals below with ideals of scaled copies of $Y_C$.
Finally, pulling back along $(\varphi_P^*)^{-1}=\varphi_{P^{-1}}^*$ and using $(\zeta^P)^{(P^{-1})_j}=\zeta_j$, we get
\[
\varphi_{P^{-1}}^*\big(\langle z_j-\zeta_j\rangle_{j\le d}\big)=\langle y^{(P^{-1})_j}-(\zeta^P)^{(P^{-1})_j}\rangle_{j\le d}=I_{C,\zeta^P}.
\]
Hence,
\[
I((Y_A\cap Y_B)\cap(\mathbb{C}^*)^n)=I_A+I_B=\bigcap_{\zeta \in Z} I_{C,\zeta^P},
\]
which concludes the proof.
\end{proof}

Combining Lemma~\ref{lemma:toricintersection} with the proof of Theorem~\ref{thm:LissajousIntersection}, we obtain the following result.

\begin{theorem}\label{thm:LAintersectiongeneral}
    Let $A\in\mathbb{Z}^{d_1\times n}$ and $B\in\mathbb{Z}^{d_2\times n}$. 
    For $J\subset[n]$, let
    $d_J=\dim(\ker A+\ker B^{-J})$, let $M^{(J)}\in\mathbb{Z}^{n\times d_J}$ constructed concatenating a $\mathbb{Z}$-basis of
    $\ker_{\mathbb{Z}}(A)+\ker_{\mathbb{Z}}(B^{-J})$ column-wise, and let $S^{(J)}=P^{(J)}M^{(J)}Q^{(J)}$ be its Smith normal form, with invariant factors $s_{J,1}, \dotsc, s_{J,{d_J}}$.
    Set
    \[
    Z_J=\{\zeta\in\C^n\, :\, \zeta_1^{s_{J,1}}=\dotsc=\zeta_{d_J}^{s_{J,d_J}}=1, \, \zeta_{k>d_J}=1\},
    \]
    and for $\zeta\in Z_J$ let $\vartheta_{J,\zeta}\in\C^n$ be such that $e^{-i{\pi\over 2} (\vartheta_{J,\zeta})_j}=(\zeta^{P^{(J)}})_j$ for all $j\in[n]$, i.e., $\vartheta_{J,\zeta}=\frac{2i}{\pi}\log(\zeta^{P^{(J)}})$. 
    Then
    \[
    {\cal L}_A\cap{\cal L}_B=\bigcup_{J\subset[n]}\ \bigcup_{\zeta\in Z_J}{\cal L}_{A\cap B^{-J},\,\vartheta_{J,\zeta}} .
    \]
\end{theorem}

\begin{example}\label{ex:intersection1213}
Let $A=\left(\begin{smallmatrix} 1&2 \end{smallmatrix}\right)$ and $B=\left(\begin{smallmatrix} 1&3 \end{smallmatrix}\right)$.
The defining equations of $\LA$ and ${\cal L}_B$ are $x_2=2x_1^2-1$ and $x_2=4x_1^3-3x_1$, respectively.
Their intersection is cut out by $4x_1^3-2x_1^2-3x_1+1=(x_1-1)(4x_1^2+2x_1-1)$, so it consists of the three points
$(1,1)$, $(\cos\tfrac{2\pi}{5},\cos\tfrac{4\pi}{5})$, and $(\cos\tfrac{4\pi}{5},\cos\tfrac{2\pi}{5})$.

Since ${\rm Row}(A)\cap{\rm Row}(B^{-J})=\{0\}$ for all $J\subset[2]$, each Lissajous variety in \Cref{thm:LAintersectiongeneral} is a single point.
For $J=\emptyset,\{1,2\}$, the lattice $M^{(J)}$ is generated by $(2,-1),(3,-1)$ and has index $1$, so $Z_J=\{{(1,1)}\}$ and we obtain the point $(1,1)$.
For $J=\{2\}$, $M^{(J)}$ is generated by $(2,-1),(3,1)$ and has lattice index $|\det\left(\begin{smallmatrix} 2 & 3\\ -1 & 1 \end{smallmatrix}\right)|=5$, with invariant factors $(1,5)$, and similarly for $J=\{1\}$. 
Hence $|Z_J|=5$.
In the torus, $Y_A\cap Y_{B^{-J}}$ is given by $y_2=y_1^2$ and $y_1^3y_2=1$, i.e., $y_1^5=1$, so the points $\zeta^{P^{(J)}}$ are $(\omega^k,\omega^{2k})$ with $\omega=e^{2\pi i/5}$ and $k=0,\dotsc,4$.
The corresponding shifts are $\vartheta_{J,\zeta}=-\frac{4k}{5}(1,2)$, and ${\cal L}_{A\cap B^{-J},\vartheta_{J,\zeta}}=\{(\cos\frac{2\pi k}{5},\cos\frac{4\pi k}{5})\}$.
The values $k=0$, $k\in\{1,4\}$ and $k\in\{2,3\}$ give exactly the three points above (see Figure~\ref{fig:intersections3.2and3.7}, on the right).
Note that \Cref{thm:LissajousIntersection}, ignoring its hypothesis, would only detect the point $(1,1)$.
\begin{figure}
    \centering
    \begin{minipage}[c]{0.45\textwidth}\centering
        \includegraphics[width=0.8\linewidth]{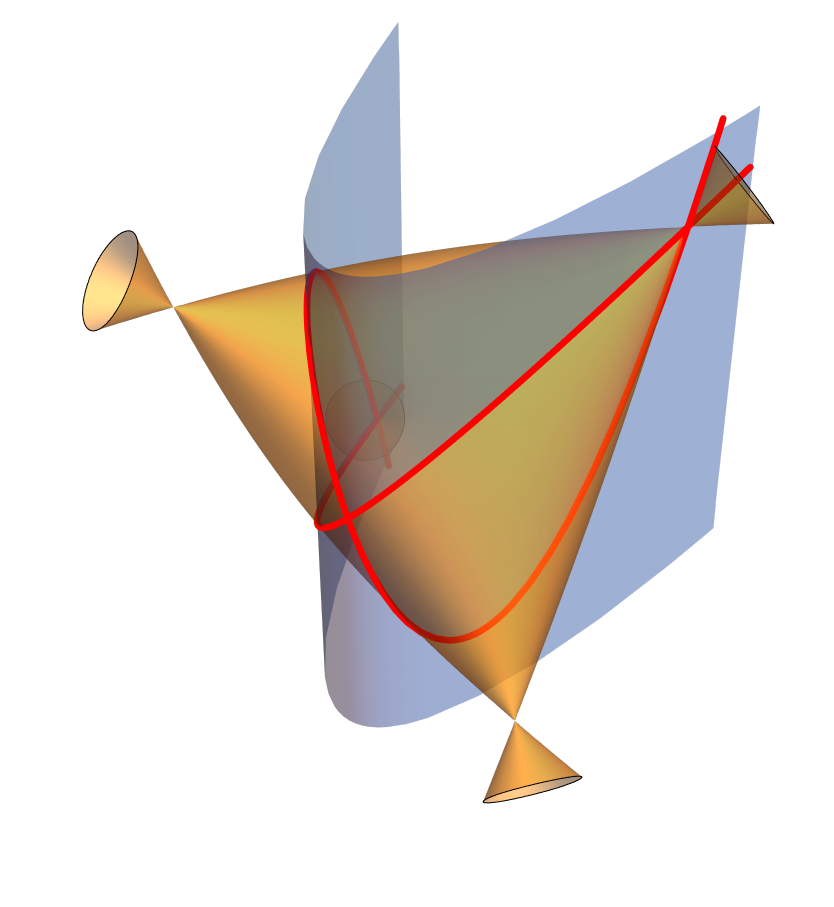}
    \end{minipage}
    \begin{minipage}[c]{0.45\textwidth}\centering
        \includegraphics[width=0.8\linewidth]{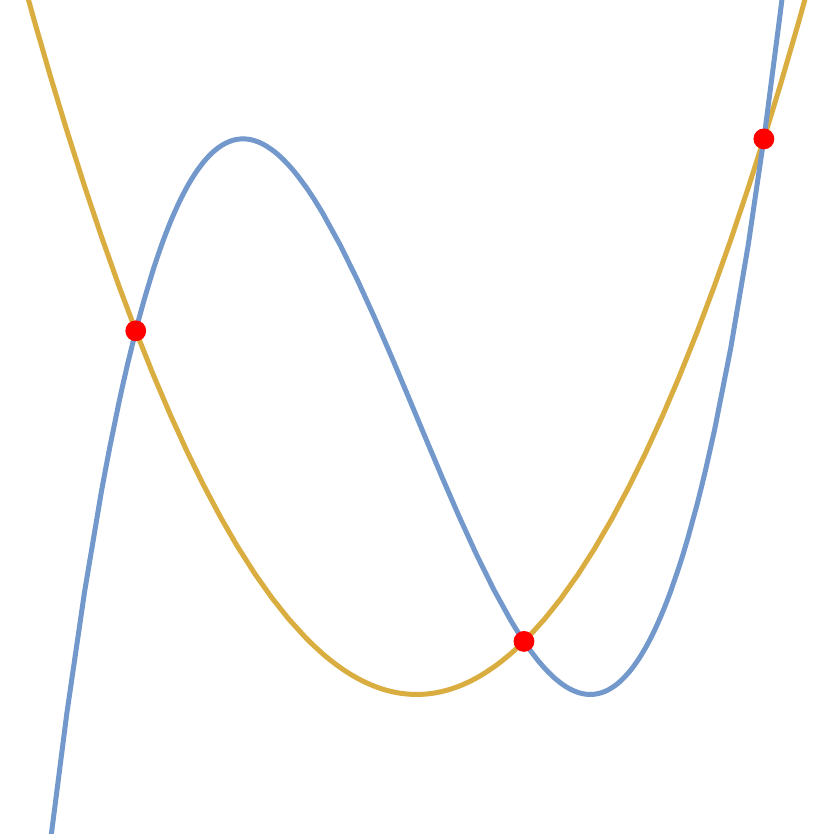}
    \end{minipage}
    \caption{On the left, the Cayley cubic (gold) and the parabolic cylinder $x_2=2x_1^2-1$ (blue) from Example~\ref{ex:intersectionlatticeindex1} intersect in two curves (red). 
    On the right, the curves $x_2=2x_1^2-1$ (gold) and $x_2=4x_1^3-3x_1$ (blue) from Example~\ref{ex:intersection1213} intersect in three points (red).}
    \label{fig:intersections3.2and3.7}
    \end{figure}
\end{example}

We finally turn to arbitrary Lissajous varieties ${\cal L}_{A_1,b_1}$ and ${\cal L}_{A_2,b_2}$.
By Theorem~\ref{thm:projY_Ab}, the shifts $b_1,b_2$ turn the toric varieties of the previous discussion into scaled toric varieties $Y_{A_1,\beta_1}$ and $Y_{A_2,\beta_2}$.
In this setting, their intersection in the torus may be empty, and its components are no longer indexed by roots of unity.
We first generalize Lemma~\ref{lemma:Yintersection} to scaled toric varieties, and then deduce the general intersection theorem, which contains Theorems~\ref{thm:LissajousIntersection} and~\ref{thm:LAintersectiongeneral} as special cases.

\begin{lemma}\label{lemma:scaledToricIntersection}
Let $K_i=\ker_{\mathbb Z}(A_i)$ ($i=1,2$) and $M=K_1+K_2$ with basis $m_1,\dots,m_d$.
Then, $(Y_{A_1,\beta_1}\cap Y_{A_2,\beta_2})\cap(\C^*)^n\neq\emptyset$ if and only if
$\beta_1^u=\beta_2^u$ for all $u\in K_1\cap K_2$.
In that case, we can define $\chi\colon M\to\C^*$, $\chi(u+w)=\beta_1^u\beta_2^w$ ($u\in K_1$, $w\in K_2$) and have
\[
I((Y_{A_1,\beta_1}\cap Y_{A_2,\beta_2})\cap(\mathbb{C}^*)^n)=\bigcap_{\zeta\in Z_c}I_{C,\zeta^P},
\]
where we identify $M=(m_1|\dotsc|m_d)\in\Z^{n\times d}$ and $S=PMQ$ is a Smith normal form, 
$Z_c=\{\zeta\in\C^n:\ \zeta_k^{s_k}=(c^Q)_k \text{ for } k\le d,\ \zeta_{k>d}=1\}$, with $c=(\chi(m_1),\dotsc,\chi(m_d))$, and $\operatorname{Row}(C)=\operatorname{Row}(A_1)\cap\operatorname{Row}(A_2)$.
\end{lemma}
\begin{proof}
    First notice that if $(Y_{A_1,\beta_1}\cap Y_{A_2,\beta_2})\cap(\C^*)^n\neq\emptyset$, then $\beta_1^u=\beta_2^u$ for all $u\in K_1\cap K_2$.
    Indeed, if $y\in(Y_{A_1,\beta_1}\cap Y_{A_2,\beta_2})\cap(\C^*)^n$ and $u\in K_1\cap K_2$, then $y^u=\beta_1^u$ and $y^u=\beta_2^u$, hence $\beta_1^u=\beta_2^u$.

    Conversely, let $\beta_{1,2}$ such that $\beta_1^u=\beta_2^u$ for all $u\in K_1\cap K_2$.
    Then, the character $\chi\colon M\to\C^*$, $\chi(v+w)=\beta_1^v\beta_2^w$ ($v\in K_1$, $w\in K_2$) is well defined.
    Indeed, for $v+w=v'+w'$, we have that $v-v'=w'-w\in K_1\cap K_2$.
    From the property of $\beta_{1,2}$, we get $\beta_1^{v-v'}=\beta_2^{v-v'}=\beta_2^{w'-w}$ and hence $\beta_1^v\beta_2^w=\beta_1^{v'}\beta_2^{w'}$, showing that $\chi(v+w)=\chi(v'+w')$.

    Let now $m=v+w\in M$, with $v\in K_1$ and $w\in K_2$. 
    Then, $y^m-\chi(m)=y^v(y^w-\beta_2^w)+\beta_2^w(y^v-\beta_1^v)\in I_{A_1,\beta_1}+I_{A_2,\beta_2}$.
    This shows that $\langle y^M - c\rangle\subseteq I_{A_1,\beta_1}+I_{A_2,\beta_2}$, where $c=(\chi(m_1),\dotsc,\chi(m_d))$.

    Conversely, let $v\in K_1$ and write $v=\sum_j\mu_jm_j$. 
    Then $y^v-\beta_1^v=\prod_j(y^{m_j})^{\mu_j}-\prod_jc_j^{\mu_j}$,
    because $\beta_1^v=\chi(v)=\prod_jc_j^{\mu_j}$. 
    Since $(y^{m_j})^{\mu_j}-c_j^{\mu_j}\in\langle y^{m_j} -c_j\rangle$ for all $\mu_j\in\mathbb Z$,
    the identity $\prod a_j-\prod b_j=\sum_k(\prod_{j<k}a_j)(a_k-b_k)(\prod_{j>k}b_j)$ implies $y^v-\beta_1^v\in\langle y^M-c\rangle$, showing that $I_{A_1,\beta_1}\subset \langle y^M-c\rangle$. 
    The same argument applies to $I_{A_2,\beta_2}$.

    Hence, 
    \[
    I_{A_1,\beta_1}+I_{A_2,\beta_2}=\langle y^M - c\rangle.
    \]
    Following the same steps as in Lemma \ref{lemma:toricintersection}, we get 
    \[
    \varphi_P^*(I_{A_1,\beta_1}+I_{A_2,\beta_2})
    =\langle (z^S)^{Q^{-1}}- c\rangle
    =\langle z^S- c^Q\rangle\subset \C[z^{\pm1}].
    \]
    Again, the ideal $\langle z^S- c^Q\rangle\subset \C[z^{\pm 1}]$ is radical, and so is $I_{A_1,\beta_1}+I_{A_2,\beta_2}$, and equals the intersection
    \[
    \langle z^S-c^Q\rangle = \bigcap_{\zeta\in Z_c} \langle z_1-\zeta_1, \dotsc, z_d-\zeta_d \rangle.
    \]
    As in the case $\beta_{1,2}=1$, from 
    \[
    \varphi_{P^{-1}}^*\big(\langle z_j-\zeta_j\rangle_{j\le d}\big)=\langle y^{(P^{-1})_j}-\zeta_j\rangle_{j\le d}=I(Y_{C,\zeta^P}\cap(\mathbb{C}^*)^n)=I_{C,\zeta^P},
    \]
    we can finally conclude
    \[
    I((Y_{A_1,\beta_1}\cap Y_{A_2,\beta_2})\cap(\mathbb{C}^*)^n)=\bigcap_{\zeta\in Z_c}I_{C,\zeta^P}.
    \]
    In particular, $Z_c\neq\emptyset$: since $\operatorname{rank}(M)=d$, all $s_k\neq0$, and as $(c^Q)_k\in\C^*$ each equation
    $\zeta_k^{s_k}=(c^Q)_k$ has $|s_k|$ solutions. For any $\zeta\in Z_c$, the point $\zeta^P\in(\C^*)^n$ lies in
    $Y_{C,\zeta^P}\cap(\C^*)^n\subseteq (Y_{A_1,\beta_1}\cap Y_{A_2,\beta_2})\cap(\C^*)^n$, which is therefore non-empty.
\end{proof}

We are now ready to state the final result, which follows from Lemma \ref{lemma:scaledToricIntersection} and only a few minor adjustments in the proof of Theorem \ref{thm:LissajousIntersection}.

\begin{theorem}\label{thm:GeneralLissajousIntersection}
Let $A_1\in\mathbb Z^{d_1\times n}$, $A_2\in\mathbb Z^{d_2\times n}$ and $b_1,b_2\in\C^n$.
For $J\subset[n]$, let $K_1=\ker_{\mathbb Z}(A_1)$, $K_2^{(J)}=\ker_{\mathbb Z}(A_2^{-J})$, and $M^{(J)}=(m_1|\dotsc|m_{d_J})$,
where $\{m_k\}_k$ is a $\mathbb Z$-basis of $K_1+K_2^{(J)}$. 
Let $S^{(J)}=P^{(J)}M^{(J)}Q^{(J)}$ be a Smith normal form with
invariant factors $s_{J,1},\dotsc,s_{J,d_J}$. 
Choose $m_k=v_k+w_k$ with $v_k\in K_1$, $w_k\in K_2^{(J)}$, and set
$\beta_{1,2}=e^{-i\frac{\pi}{2}b_{1,2}}$, $c_J=(\beta_1^{v_k}(\beta_2^{-J})^{w_k})_{k\le d_J}$ and
\[
Z_J=\{\zeta\in\C^n:\ \zeta_k^{s_{J,k}}=(c_J^{Q^{(J)}})_k\ \text{for } k\le d_J,\ \zeta_{k>d_J}=1\}.
\]
Say that $J$ satisfies (P) if $b_1-b_{2,-J}\in\Row(A_1)+\Row(A_2^{-J})+4\Z^n$.
Let $C^{(J)}$ be such that $\Row(C^{(J)})=\Row(A_1)\cap\Row(A_2^{-J})$, and for $\zeta\in Z_J$ let
$\vartheta_{J,\zeta}=\frac{2i}{\pi}\log(\zeta^{P^{(J)}})$. 
Then
\[
{\cal L}_{A_1,b_1}\cap{\cal L}_{A_2,b_2}=\bigcup_{\substack{J\subset[n]\\ J\text{ s.t. (P)}}}\ \bigcup_{\zeta\in Z_J}{\cal L}_{C^{(J)},\vartheta_{J,\zeta}}.
\]
\end{theorem}
\begin{proof}
    Following the proof of Theorem \ref{thm:LissajousIntersection}, combined with the fact that $z\in Y_{A_2,\beta_2}$ if and only if $y\in Y_{A_2^{-J},\beta_2^{-J}}$, we get 
    \[
    {\cal L}_{A_1,b_1}\cap{\cal L}_{A_2,b_2}=\bigcup_{J\subset[n]}\pi_x\big({\cal Y}_{A_1,b_1}\cap {\cal Y}_{A_2^{-J},b_{2,-J}}\big).
    \]
    Let $J\subset[n]$ and apply Lemma \ref{lemma:scaledToricIntersection} to data $(A_1,b_1)$ and $(A_2^{-J},b_{2,-J})$: the intersection $Y_{A_1,\beta_1}\cap Y_{A_2^{-J},\beta_2^{-J}}\cap (\C^*)^n\ne \emptyset$ precisely when $\beta_1^u=(\beta_2^{-J})^u$ for all $u\in K_1\cap K_{2}^{(J)}$. 
    That is, $e^{-i{\pi\over 2} (b_1-b_{2,-J})\cdot u}={\bf 1}$, and so $(b_1-b_{2,-J})\cdot u\in 4\Z$ for all $u\in K_1\cap K_2^{(J)}=\Z^n\cap (\Row(A_1)+\Row(A_2^{-J}))^\perp$.
    The lattice $K_1\cap K_2^{(J)}$ is saturated, so the condition is equivalent to (P).
    Hence, the intersection is non-empty if and only if (P) holds for $J$, and in that case,
    \[
    I((Y_{A_1,\beta_1}\cap Y_{A_2^{-J},\beta_2^{-J}})\cap(\mathbb{C}^*)^n)=\bigcap_{\zeta\in Z_J}I_{C^{(J)},\zeta^{P^{(J)}}}.
    \]
    By definition of ${\cal Y}_{A,b}$, we have
    $\mathcal Y_{A_1,b_1}\cap\mathcal Y_{A_2^{-J},b_{2,-J}}=\bigcup_{\zeta\in Z_J}\mathcal Y_{C^{(J)},\vartheta_{J,\zeta}}$.
    Projecting onto $\C^n_x$ and taking the union over $J\subset [n]$ gives the proof. 
\end{proof}

\begin{example}\label{ex:intersectionCASA}
Let $A = \left ( \begin{smallmatrix} 1 & 0 & -1 \\ -1 & 1 & 0 \end{smallmatrix} \right)$, the reduced incidence matrix of the triangle graph $C_3$, so that $\ker(A) = \langle u\rangle$ with $u=(1,1,1)$.
The surface $\LA={\cal L}_{A,{\bf 0}}$ is the Cayley cubic, which is the algebraic boundary of the elliptope ${\cal E}_3$, and the surface ${\cal L}_{A,{\bf 1}}$ plays a fundamental role in capturing the equilibria of the Kuramoto model of the triangle graph.
The defining polynomial of ${\cal L}_{A}$ is the $3$rd cycle polynomial $\Gamma_3=x_1^2+x_2^2+x_3^2-2x_1x_2x_3-1$, and that of ${\cal L}_{A,{\bf 1}}$ is the sextic
$x_1^4+x_2^4+x_3^4-2\sum_{i<j}x_i^2x_j^2+4x_1^2x_2^2x_3^2$.
We apply Theorem~\ref{thm:GeneralLissajousIntersection} with $A_1=A_2=A$, $b_1={\bf 0}$, $b_2=\bf 1$, so that $\beta_1={\bf 1}$ and $\beta_2=(-i,-i,-i)$.

For $J=\emptyset,[3]$ we have $K_1\cap K_2^{(J)}=\langle u\rangle$ and $(b_1-b_{2,-J})\cdot u=\mp3\notin4\Z$, so (P) fails.
Note that for $J=\emptyset$ the lattice $\ker_{\Z}(A)+\ker_{\Z}(A)$ has index $1$, and yet it contributes nothing because of the scalings: $y^u=1$ and $y^u=\beta_2^u=i$ are incompatible.
For $|J|\in\{1,2\}$, the vectors $u$ and $u_{-J}$ are independent, so $K_1\cap K_2^{(J)}=0$ and (P) holds trivially.
Moreover, $\ker_{\Z}(A)+\ker_{\Z}(A^{-J})$ has index $2$ (for $J=\{1\}$, the $2\times2$ minors of $(u\,|\,u_{-J})$ are $2,2,0$), so $Z_J$ has two elements.
For $J=\{1\}$ one finds $C^{(J)}=(0\ 1\ {-1})$, and the two components are given, for instance, by $\vartheta=(\tfrac12,-\tfrac14,-\tfrac14)$ and $\vartheta=(-\tfrac32,\tfrac34,\tfrac34)$.
The sets $J$ and $[3]\smallsetminus J$ give the same components (they have the same $x$-image), and $J=\{2\},\{3\}$ give the same curves with $x_2,x_3$ in place of $x_1$.
Altogether, ${\cal L}_{A}\cap{\cal L}_{A,{\bf 1}}$ is the union of six ellipses (see Figure \ref{fig:intersection}),
\[
{\cal L}_{A}\cap{\cal L}_{A,{\bf 1}}=\bigcup_{j=1}^{3}\ \bigcup_{\varepsilon=\pm1}\Big\{x_j=\varepsilon\tfrac{\sqrt2}{2},\ x_k^2+x_l^2-\varepsilon\sqrt2\,x_kx_l=\tfrac12\Big\},\qquad\{j,k,l\}=\{1,2,3\}.
\]
\begin{figure}
    \centering
    \includegraphics[height = 8cm]{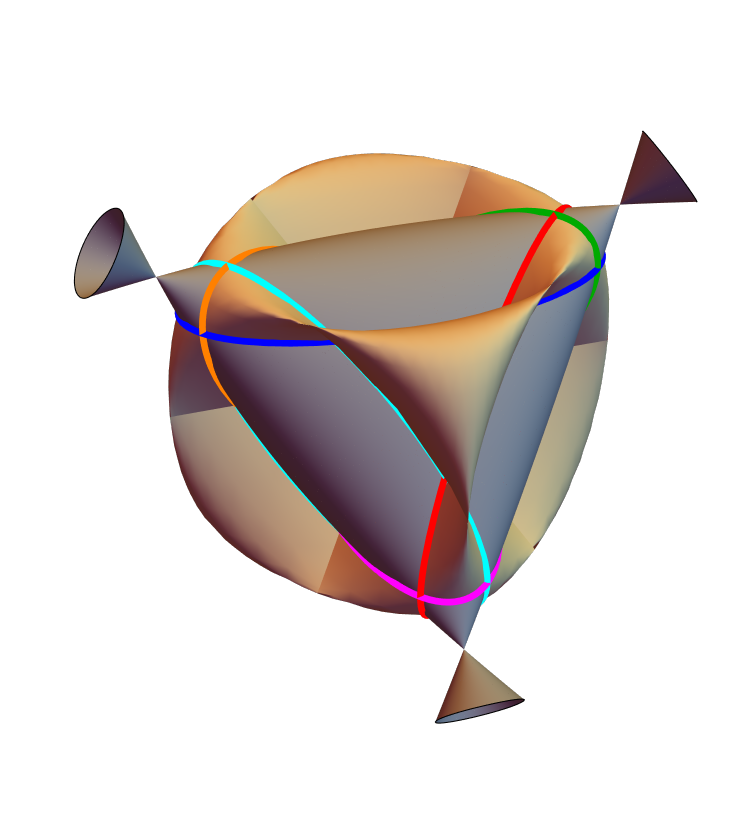} 
    \includegraphics[height = 8cm]{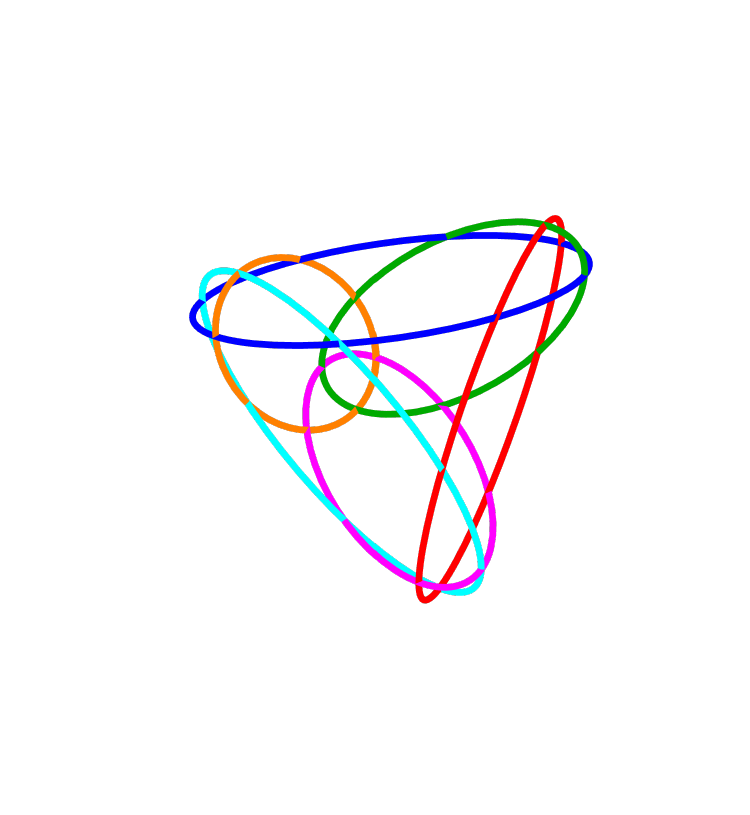}
    \caption{The intersection ${\cal L}_{A}\cap{\cal L}_{A,{\bf 1}}$, with $A=A(C_3)$, is a union of six ellipses (Example \ref{ex:intersectionCASA}).}
    \label{fig:intersection}
    \end{figure}
\end{example}

\section{Lissajous hypersurfaces}\label{sec:hypersurfaces}
If the matrix $A\in\Z^{(n-1)\times n}$ is of rank $n-1$, then $L_{A,b}$ is a hyperplane and by Theorem \ref{thm:mainintro} the Lissajous variety $\LAb=\cos(L_{A,b})$ is a hypersurface.
Our goal in this section is to describe the defining polynomial of Lissajous hypersurfaces $\LAb$, when $b\in\C^n$ is generic.

A set-theoretical determinantal representation of Lissajous hypersurfaces is derived in \cite{mascarin2025lissajousvarieties} from the characterization provided by Theorem \ref{thm:projY_Ab} (see \cite[Theorem 3.1]{mascarin2025lissajousvarieties}).
The construction proceeds as follows. 
Let $v$ be a primitive generator of $\ker(A)\cap\Z^n$, so that $I_{A,\beta}=I(Y_{A,\beta}\cap(\C^*)^n )$ is generated by $g_{A,b}=y^v-\beta^v$, where $\beta=\exp(-i\frac{\pi}{2}b)$ as usual.
Let $K = \overline{\mathbb{C}(x)}$ be the algebraic closure of the field of rational functions in $x = (x_1, \ldots, x_n)$, and consider the ideal: 
\[ 
J \, = \, \langle 2 \, x_1 - y_1 - y_1^{-1}, \, \ldots, \, 2\, x_n - y_n - y_n^{-1}\rangle \, \, \subset \, \, K[y^{\pm}] = K[y_1^{\pm 1}, \ldots, y_n^{\pm 1}]. 
\]
The affine variety $V(J)$ defined by $J$ in $(K^*)^n$ is zero-dimensional and consists of $2^n$ distinct points,
\[
V(J)\, = \, \left\{\, \big(x_j+\epsilon_j\sqrt{x_j^2-1}\big)_{j=1}^n\, :\, \epsilon\in\{\pm1\}^n\, \right\}.
\]
Consider the algebra 
\[
{\cal A} \, = \, \frac{K[y^{\pm 1}]}{\langle f_1,\dotsc,f_n\rangle}\, \simeq \, \frac{K[y_1^{\pm 1}]}{\langle f_1\rangle}\otimes_K\, \cdots \, \otimes_K\frac{K[y_n^{\pm 1}]}{\langle f_n\rangle} \, , 
\]
where $f_j=2  x_j - y_j - y_j^{-1}$.
A $K$-linear endomorphism $M_g : {\cal A} \rightarrow {\cal A}$ is obtained by considering the multiplication with an element $g \in K[y^{\pm 1}]$, where $M_g([h]) = [gh]$ (with $[h]$ the residual class of $h$), which is represented by a $2^n \times 2^n$ matrix with entries in $K$, once a basis is fixed.
It follows from the eigenvalue method for polynomial root finding (see as a reference \cite[Chapter 4, Proposition 2.7]{CoxLittleOShea2005UsingAG})  that the eigenvalues of $M_g$ are the evaluations $g(\xi)$ of $g$ at the roots $\xi\in V(J)$. 
In particular, we have
\begin{align}\label{eq:detMg}
\det (M_g)=  \prod_{\xi\in V(J)} g(\xi).
\end{align}
By Theorem \ref{thm:projY_Ab}, $x\in\ML_{A,b}$ if and only if there exists $y\in V(J)$ such that $y\in Y_{A,\beta}\cap (\C^*)^n=V(g_{A,b})$; equivalently, $x\in\ML_{A,b}$ if and only if $\prod_{y \in V(J)}g_{A,b}(y)=0$. 
Therefore, using (\ref{eq:detMg}) for the polynomial $g = g_{A,b}$, we obtain the following {\em determinantal representation} of $\ML_{A,b}$,
\begin{align}\label{eq:detrep}
\ML_{A,b}=\{x\in\C^n: \det(M_{g_{A,b}})=0\}.
\end{align}
Denote with $\Gamma_{A,b}$ the irreducible generator of the prime ideal of ${\cal L}_{A,b}$ (up to a scalar).
By \cite[Proposition 3.5]{mascarin2025lissajousvarieties}, we know that $\det(M_{g_{A,b}}) = c \,  \Gamma_{A,b}^{\deg \pi_{A,b}}$ for some~$c \in \mathbb{C}^*$, with $\pi_{A,b}$ as in Theorem \ref{thm:projY_Ab}.
The above discussion motivates the following definition.
\begin{definition}
    Let $(A,b)$ as before and $v$ be the primitive generator of $\ker(A)\cap\Z^n$ with first non-zero entry positive. 
    We define the polynomial $h_{A,b}$ to be the determinant of the multiplication matrix $M_{g_{A,b}}$ from the determinantal representation (\ref{eq:detrep}), namely
    \begin{equation}\label{eq:h_Ab}
        h_{A,b} = \prod_{\xi\in V(J)} g_{A,b}(\xi) = \prod_{\epsilon\in\{\pm1\}^n}\left(\prod_{j=1}^n\left(x_j+\epsilon_j\sqrt{x_j^2-1}\right)^{v_j}-\beta^v\right)
    \end{equation}
\end{definition}
As an immediate consequence of \cite[Theorem 2.4]{mascarin2025lissajousvarieties} we get the following result.
\begin{proposition}
    If $\Z A=\Z^{n-1}$ and $v \cdot b$ is not an even integer, then $\deg(\pi_{A,b})=2^{{\rm CL}_A}$, where ${\rm CL}_A$ is the number of coloops of $A$ (equivalently, the number of zero entries of $v$). 
    Hence, $h_{A,b}=c\Gamma_{A,b}^{2^{{\rm CL}_A}}$, with $\Gamma_{A,b}$ the irreducible generator of $I(\LAb)$ and $c\in\C^*$.
\end{proposition}
In particular, if $v$ has no zero entry, then $h_{A,b}$ (for $b\in\C^n$ generic) is indeed the irreducible generator of the ideal of $\LAb$.
This motivates studying its properties.
First, note from (\ref{eq:LissDeg}) that, if $A$ has no lattice index, the degree of $h_{A,b}$ is given by $\deg(h_{A,b})=(n-1)!{\rm vol}(P_A)$.
In this section we derive more refined information, namely its univariate degree and leading coefficient, together with its roots as a univariate polynomial.
Before proceeding, we introduce some notation to simplify the discussion that follows.
As in \cite[Section 3.1]{Laurent2026AlgebraicBoundary}, we define the following algebraic functions of $x$:
\begin{align}\label{eq:zepsilon}
z_\epsilon\, =\, \prod_{j=1}^n \left( x_j + \epsilon_j \sqrt{x_j^2-1} \right)^{|\epsilon_j|} \ \text{ for } \epsilon \in \{\pm1, 0\}^n.
\end{align}
Note that $z_{-\epsilon} = z_\epsilon^{-1}$ and $z_{\bf 0}=1$.
Moreover, if $\delta$, $\delta'$ are sequences whose concatenation $(\delta,\delta')$ lies in $\{\pm 1,0\}^n$, then $z_{\delta,\delta'}=z_{\delta,\azero }\cdot z_{\azero, \delta'}$ holds. 
For $v\in\Z^n$, we write
\begin{equation}\label{eq:zepsilonv}
z_{\epsilon}^v\, =\, \prod_{j=1}^n \left( x_j + \epsilon_j \sqrt{x_j^2-1} \right)^{|\epsilon_j|v_j}.
\end{equation}
We stress that $z_\epsilon^v$ is not a power of the scalar $z_\epsilon$, whereas for
vectors such as $\beta\in(\C^*)^n$ we keep the usual multi-index power
$\beta^v=\beta_1^{v_1}\cdots\beta_n^{v_n}$. Both conventions are used throughout this section.
In particular, $z_{-\epsilon}^v=(z_\epsilon^v)^{-1}$, $z_{\delta,\delta'}^v=z_{\delta,\azero}^v\, z_{\azero,\delta'}^v$,
and if $\epsilon$ is supported only in the $j$-th coordinate, then $z_\epsilon^v=z_\epsilon^{v_j}$.
Using \eqref{eq:zepsilonv}, we can rewrite $h_{A,b}$ from \eqref{eq:h_Ab} as
\begin{equation}\label{eq:h_Ab2}
    h_{A,b}=\prod_{\epsilon\in \{\pm1\}^n}(z_{\epsilon}^v-\beta^v).
\end{equation}
With this notation, we are now ready to state and prove the main result of this section, which generalizes \cite[Proposition 3.6]{Laurent2026AlgebraicBoundary}.

\begin{theorem}\label{thm:LissHypersurfaces}
Let $A\in\Z^{(n-1)\times n}$ have rank $n-1$, let $b\in\C^n$, and let $v$ be the primitive
generator of $\ker(A)\cap\Z^n$ whose first non-zero entry is positive.
Fix $j\in[n]$ with $v_j\neq 0$, and regard $h_{A,b}$ as a polynomial in $x_j$ with
coefficients in $\C[x_l : l\neq j]$.
Then $h_{A,b}$ has degree $2^{n-1}|v_j|$ in $x_j$, and its leading coefficient is the
non-zero constant $(2^{|v_j|}\beta^v)^{2^{n-1}}$.

Moreover, for $j=1$, the $2^{n-1}|v_1|$ roots of $h_{A,b}$ in $x_1$ (counted with multiplicity)
in the algebraic closure of $\C(x_2,\dotsc,x_n)$ are
\[
r^{(k,\delta)}=\tfrac12\big(\xi^{k}f_\delta+\xi^{-k}f_\delta^{-1}\big),
\qquad \delta\in\{\pm1\}^{n-1},\ k=0,\dotsc,|v_1|-1,
\]
where $\xi=e^{2\pi i/|v_1|}$ is a primitive $|v_1|$-th root of unity and $f_\delta$ is any
fixed solution of $f_\delta^{|v_1|}=\beta^{-v}z_{0\delta}^{v}$.
Equivalently, for each $\delta\in\{\pm1\}^{n-1}$ the roots $r^{(k,\delta)}$,
$k=0,\dotsc,|v_1|-1$, are the solutions in $x_1$ of
\begin{equation}\label{eq:ChebyshevEquation}
2\,T_{|v_1|}(x_1)=\beta^{-v}z_{0\delta}^{v}+\beta^{v}z_{0\delta}^{-v},
\end{equation}
where $T_{|v_1|}$ is the $|v_1|$-th Chebyshev polynomial of the first kind.

For general $j$ the same holds, with the zero in $z_{0\delta}$ placed in position $j$.
\end{theorem}

Before proceeding with the proof, we need the following short technical lemma. 
For convenience we index $z_\epsilon$ by $\{\pm,0\}$ instead of $\{\pm 1,0\}$.

\begin{lemma}\label{lemma:ChebyshevPoly}
    With the above notation,
    \[
    z_{+{\bf 0}}^v+z_{-{\bf 0}}^v=2T_{|v_1|}\in\R[x],
    \]
    where $T_{|v_1|}$ denotes the $|v_1|$-th Chebyshev polynomial of the first kind.
\end{lemma}
\begin{proof}
    Observe that $z_{\pm{\bf 0}}=x_1\pm\sqrt{x_1^2-1}$ are the two solutions to the equation $z^2-2x_1z+1=0$. 
    In particular, for $v_1>0$,
    \[
    z_{\pm{\bf 0}}^{v_1+1}=z_{\pm{\bf 0}}^{v_1-1}z_{\pm{\bf 0}}^2 = 2x_1z_{\pm{\bf 0}}^{v_1}-z_{\pm{\bf 0}}^{v_1-1}.              
    \]
    Setting $Z_{v_1}=\frac{z_{+{\bf 0}}^v + z_{-{\bf 0}}^v  }{2}$, the recursion $Z_{v_1+1}=2x_1Z_{v_1}-Z_{v_1-1}$ holds.
    In particular, $Z_0=1$ and $Z_1=x_1$, and hence by induction $Z_{v_1}(x_1)=T_{v_1}(x_1)$.
    The case $v_1<0$ is analogous, as $z_{\pm{\bf 0}}^{-1}=z_{\mp{\bf 0}}$. 
\end{proof}
\begin{proof}[Proof of Theorem \ref{thm:LissHypersurfaces}]
    We derive the result for $j=1$; for the remaining variables the procedure is identical.
    Without loss of generality, we assume $v_1>0$.
    We manipulate $h_{A,b}$ from the compact expression (\ref{eq:h_Ab2}).
    \begin{align*}
        h_{A,b} &=\, \prod_{\epsilon\in \{\pm\}^n}(z_{\epsilon}^v-\beta^v) = \prod_{\delta\in \{\pm\}^{n-1}}(z_{+{\bf 0}}^vz_{0\delta}^v-\beta^v)(z_{-{\bf 0}}^vz_{0\delta}^v-\beta^v)\,=\\
        &=\, \prod_{\delta\in \{\pm\}^{n-1}}(z_{0\delta}^{2v}-\beta^vz_{0\delta}^v(z_{+{\bf 0}}^v+z_{-{\bf 0}}^v)+\beta^{2v})\,
        =\, \prod_{\delta\in \{\pm\}^{n-1}}(z_{0\delta}^{2v}-2\beta^vT_{v_1}z_{0\delta}^v+\beta^{2v})
    \end{align*}
    Here we used the property $z_{\epsilon}^v=z_{\epsilon_1{\bf 0}}^vz_{0\delta}^v$ and Lemma \ref{lemma:ChebyshevPoly}.
    To simplify notation, we write for the moment $p_{v_1}=2T_{v_1}$ and $s_\delta^v=\beta^{-v}z_{0\delta}^v+\beta^vz_{0\delta}^{-v}$.
    Factoring out $\beta^vz_{0\delta}^v$ from each term of the product, we get:
    \begin{align*}
        h_{A,b} \,&=\, \prod_{\delta\in \{\pm\}^{n-1}}\beta^vz_{0\delta}^v(s_\delta^v-p_{v_1})\,
        =\,\prod_{\delta\in \Delta}\beta^vz_{0,\delta}^v(s_\delta^v-p_{v_1})\beta^vz_{0,-\delta}^v(s_{-\delta}^v-p_{v_1})\,=\\
        &=\,\prod_{\delta\in \Delta}\beta^{2v}(s_\delta^v-p_{v_1})(s_{-\delta}^v-p_{v_1})\,
        =\, (\beta^{2v})^{2^{n-2}}\prod_{\delta\in \Delta}\big(p_{v_1}^2-p_{v_1}(s_{\delta}^v+s_{-\delta}^v)+c_\delta\big)
    \end{align*}
    where we partition $\{\pm1\}^{n-1}=\Delta\cup(-\Delta)$, with $\Delta=\{\delta\in\{\pm1\}^{n-1}\,:\,\delta_1=+1\}$, and write $c_\delta=z_{0,-\delta}^{2v}+z_{0,-\delta}^{-2v}+\beta^{2v}+\beta^{-2v}$, which does not depend on $x_1$.
    Since $p_{v_1}=2T_{v_1}$, and Chebyshev polynomials of the first kind have degree $v_1$ and leading coefficient $2^{v_1-1}$, we get:
    \[
    h_{A,b}=(\beta^{2v})^{2^{n-2}}\prod_{\delta\in \Delta}\big(2^{2v_1}x_1^{2v_1}\, +\, {\rm l.o.t.}\big)=(2^{2^{n-1}v_1}\beta^{2^{n-1}v})\cdot x_1^{2^{n-1}v_1}\, +\, {\rm l.o.t.},
    \]
    where l.o.t. denotes lower order terms in $x_1$.
    This proves the first part of the claim.

    Since $h_{A,b} = \prod_{\delta\in \{\pm\}^{n-1}}\beta^vz_{0\delta}^v(s_\delta^v-p_{v_1})$, the second part of the statement follows by solving the rational equation $s_\delta^v-p_{v_1}=0$, which gives equation (\ref{eq:ChebyshevEquation}).
    From Lemma \ref{lemma:ChebyshevPoly}, we know that $2T_{v_1}=z_{+{\bf 0}}^v+z_{+{\bf 0}}^{-v}$ (recall that $z_\epsilon^{-1}=z_{-\epsilon})$.
    If we write $y_1=z_{+{\bf 0}}^v$, then from $2T_{v_1}=s_\delta^v$ we get
    \[
    y_1^2-s_\delta^vy_1+1=0,
    \]
    and hence $y_1=\beta^{\mp v}z_{0\delta}^{\pm v}$.
    Since $x_1=\frac{1}{2}(z_{+{\bf 0}}+z_{+{\bf 0}}^{-1})$, we finally get
    \[
    x_1=\frac{1}{2}\big( \xi_{v_1}^k (y_1)^{1\over v_1} +\xi_{v_1}^{-k} (y_1)^{-{1\over v_1}} \big)=\frac{1}{2}\big( \xi_{v_1}^k (\beta^{-v}z_{0\delta}^v)^{1\over v_1} +\xi_{v_1}^{-k} (\beta^{-v} z_{0\delta}^v)^{-{1\over v_1}} \big),
    \]
    where $\xi_{v_1}^k=e^{2k\pi i\over v_1}$, for $k=0,\dotsc,v_1-1$ denote the $v_1$-th roots of unity. 
    Here $(\beta^{-v}z_{0\delta}^v)^{1\over v_1}$ denotes any algebraic function $f$ satisfying the equation $f^{v_1}=\beta^{-v}z_{0\delta}^v$.
\end{proof}

\begin{example}\label{ex:hypersurface}
    The matrix $A=\left ( \begin{smallmatrix} 3&-1&0\\5&0&-1\end{smallmatrix} \right)$ from Example \ref{ex:intro} and $b$ generic gives $h_{A,b}=\Gamma_{A,b}$ the irreducible generator of $I(\LAb)$ (up to a constant).
    Here $v=(1,3,5)$, so the polynomial is obtained as $h_A=\prod_{\epsilon\in\{\pm1\}^3}(z_{\epsilon_1}\,z_{\epsilon_2}^{3}\,z_{\epsilon_3}^{5}-\beta_1\beta_2^3\beta_3^5)$.
    By Theorem \ref{thm:LissHypersurfaces}, $h_{A,b}$ has $x_j$-degrees $(4, 12, 20)$, and corresponding leading coefficients $((2\beta^v)^4, (2^3\beta^v)^4, (2^5\beta^v)^4)$.
    When viewed as a univariate polynomial in, say, $x_3$, it reads:
    \[
    h_{A,b} = (2^5\beta^v)^4\prod_{\delta\in\{\pm 1\}^2} \prod_{k=0}^4 \left(x_3-\frac{1}{2}\Big( e^{\frac{2}{5}k\pi i}(\beta^{-v}z_{\delta0}^v)^{\frac{1}{5}} + e^{-\frac{2}{5}k\pi i}(\beta^{-v}z_{\delta0}^v)^{-\frac{1}{5}}\Big)\right)
    \]
    The hypersurface ${\cal L}_{A,{\bf 1}}$ is shown in Figure \ref{fig:Sin135}. 
    Here $b=\bf 1$ satisfies the genericity condition so $h_{A,{\bf 1}}$ is its defining polynomial (up to a constant).
    \begin{figure}
    \centering
    \includegraphics[height = 6cm]{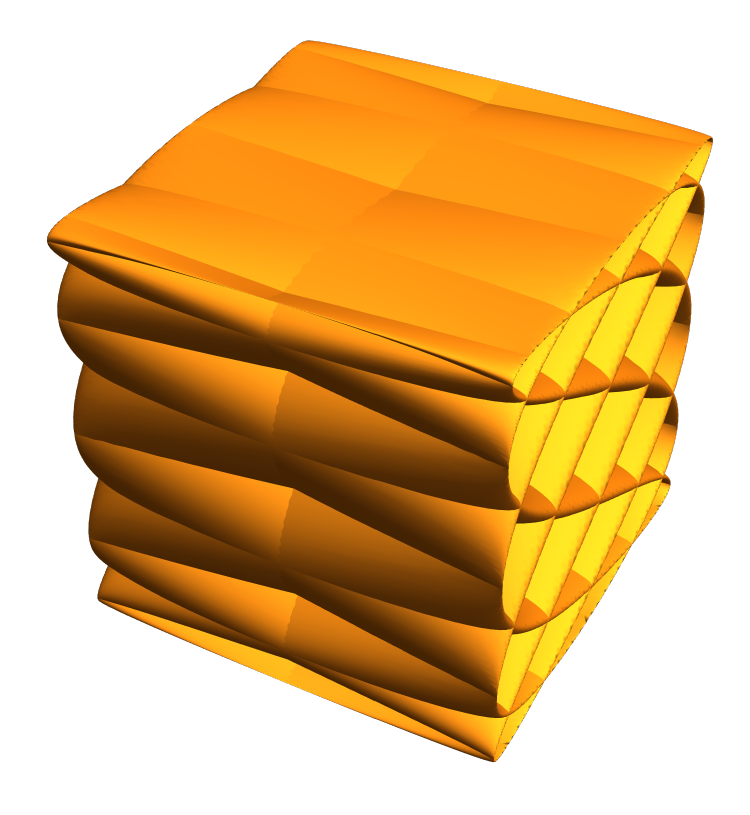} 
    \quad \quad
    \includegraphics[height = 6cm]{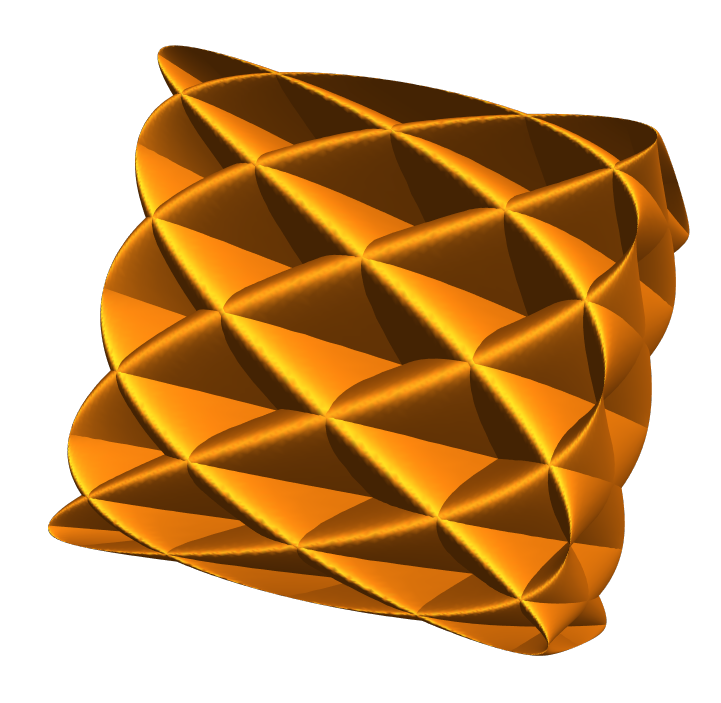}
    \caption{Two views of the Lissajous hypersurface ${\cal L}_{A,\bf 1}$ from Example \ref{ex:hypersurface}.}
    \label{fig:Sin135}
    \end{figure}
\end{example}

Let $A_1,A_2\in\Z^{(n-1)\times n}$ have rank $n-1$ and $b_1,b_2\in\C^n$.
Let $\Gamma_{A_1,b_1},\Gamma_{A_2,b_2}\in\C[x_1,\dotsc,x_n]$ be
the defining polynomials of ${\cal L}_{A_1,b_1}$ and ${\cal L}_{A_2,b_2}$. 
For $b_i$ generic and $A_i$ coloops-free we have already observed that these are precisely the polynomials $h_{A_i,b_i}$ (for $i=1,2$).
Suppose that the variable $x_n$ appears in both.
By Theorem~\ref{thm:LissHypersurfaces}, the leading coefficient of $\Gamma_{A_i,b_i}$ with respect to $x_n$ is a non-zero constant $\lambda_i$.
Writing $\Gamma_{A_1,b_1}=\lambda_1\prod_j(x_n-\eta_j)$ and $\Gamma_{A_2,b_2}=\lambda_2\prod_k(x_n-\zeta_k)$ over the algebraic closure of $\C(x_1,\dotsc,x_{n-1})$, by the Poisson formula for the resultant (see e.g., \cite[Ch.~3, Eq.~(1.4)]{CoxLittleOShea2005UsingAG}) we get
\[
\operatorname{Res}_{x_n}(\Gamma_{A_1,b_1},\Gamma_{A_2,b_2})=\lambda_1^{d_2}\lambda_2^{d_1}\prod_{j,k}(\eta_j-\zeta_k),\quad d_i=\deg_{x_n}(\Gamma_{A_i,b_i}).
\]
Since the leading coefficients never vanish, for every $\tilde{x}^*\in\C^{n-1}$, $\operatorname{Res}_{x_n}(\Gamma_{A_1,b_1},\Gamma_{A_2,b_2})(\tilde{x}^*)=0$ if and only if $\Gamma_{A_1}(\tilde{x}^*,x_n)$ and $\Gamma_{A_2}(\tilde{x}^*,x_n)$ have a common root, with $\tilde{x}=(x_1,\dotsc,x_{n-1})$. 
Hence
\[
V\big(\operatorname{Res}_{x_n}(\Gamma_{A_1,b_1},\Gamma_{A_2,b_2})\big)=\pi_{x_n}\big({\cal L}_{A_1,b_1}\cap{\cal L}_{A_2,b_2}\big),
\]
where $\pi_{x_n}\colon\C^n\to\C^{n-1}$ forgets $x_n$.
For more details on resultants and their properties, see \cite[Ch.~3, Sec.~6]{cox2025ideals}, and \cite[Ch.~3]{CoxLittleOShea2005UsingAG}.
If ${\cal L}_{A_1,b_1}\neq{\cal L}_{A_2,b_2}$, then $\Gamma_{A_1,b_1},\Gamma_{A_2,b_2}$ are distinct irreducible polynomials by definition, so they have no common factor of positive degree in $x_n$ and $\operatorname{Res}_{x_n}(\Gamma_{A_1,b_1},\Gamma_{A_2,b_2})\not\equiv0$.

By Theorem~\ref{thm:GeneralLissajousIntersection} of Section~\ref{sec:intersections}, ${\cal L}_{A_1,b_1}\cap{\cal L}_{A_2,b_2}$ is a finite union of Lissajous varieties ${\cal L}_{C^{(J)},\vartheta_{J,\zeta}}$.
Since $\cos$ acts coordinate-wise, $\pi_{x_n}({\cal L}_{C,\vartheta})=\cos(\pi_{x_n}({\rm Row}(C)-\tfrac{\pi}{2}\vartheta))$ is again a Lissajous variety, and it is irreducible.
Therefore
\[
V\big(\operatorname{Res}_{x_n}(\Gamma_{A_1,b_1},\Gamma_{A_2,b_2})\big)\,=\,\bigcup_{J,\zeta}\pi_{x_n}\big({\cal L}_{C^{(J)},\vartheta_{J,\zeta}}\big),
\]
and the resultant is in general reducible. 
Its irreducible factors are, up to constants, the defining polynomials of those $\pi_{x_n}({\cal L}_{C^{(J)},\vartheta_{J,\zeta}})$ that are hypersurfaces of $\C^{n-1}$ (the components of higher codimension are contained in these, because $V(\operatorname{Res})$ has codimension one), possibly with multiplicity. 
Different pairs $(J,\zeta)$ may give the same factor.

In the setting of Corollary~\ref{cor:intersection} with $b_1=b_2=\bf 0$, the resultant detects the single irreducible factor defining $\pi_{x_n}({\cal L}_{A_1\cap A_2})$ (up to multiplicity).
This construction gives a resultant representation for all Lissajous hypersurfaces of type ${\cal L}_A$, that can be generalized to any $\LAb$, under stronger hypotheses that guarantee a single component in the intersection.

To get a resultant representation for a Lissajous hypersurface ${\cal L}_A\subset\C^{n-1}$, one can proceed as follows.
Let $w$ generate $\ker_{\Z}(A)$, split $\operatorname{supp}(w)=S_1\sqcup S_2$,
and let $w_i$ be the restriction of $w$ to $S_i$. 
Choose $A_1,A_2$ with $\ker_{\Z}(A_1)=\langle(w_1,-1)\rangle$ and $\ker_{\Z}(A_2)=\langle(w_2,1)\rangle$.
Here each $A_i$ is $|S_i|\times(|S_i|+1)$ and acts on the variables $S_i\cup \{n\}$. 
Regard ${\cal L}_{A_i}$ ($i=1,2$) as cylinders in $\C^n$, free in some of the coordinates.
Their defining equations share only $x_n$, the lattice index of $\ker_{\Z}(A_1)+\ker_{\Z}(A_2)$ equals $\gcd(w)=1$, and $\pi_{x_n}({\cal L}_{A_1\cap A_2})={\cal L}_A$.
Hence ${\cal L}_A$ is cut out by $\operatorname{Res}_{x_n}(\Gamma_{A_1},\Gamma_{A_2})$, where $\Gamma_{A_1}$ involves only $x_n$ and the variables in $S_1$, and $\Gamma_{A_2}$ only involves $x_n$ and the variables in $S_2$.

For $b_i=\bf0$ the polynomial $h_{A_i}=h_{A_i,\bf 0}$ is the square of the reduced equation, $h_{A_i}=c_i\Gamma_{A_i}^2$ (see the proof of \cite[Proposition 2.6]{mascarin2025lissajousvarieties}), so to avoid a power factor $4$, the resultant must be computed from the polynomials $\Gamma_{A_i}$: indeed, $\operatorname{Res}_{x_n}(h_{A_1},h_{A_2})=\operatorname{Res}_{x_n}(\Gamma_{A_1},\Gamma_{A_2})^4$ (up to a constant).
This construction in the special case where $A=A_n$ is the reduced incidence matrix of the cycle graph $C_n$ is exploited in detail, showing that the resultant is irreducible, in \cite[Section 3.2]{Laurent2026AlgebraicBoundary} to obtain a resultant representation of the \emph{cycle polynomial}.
More generally, every Lissajous hypersurface is obtained from a cycle polynomial by a Chebyshev substitution.

\begin{proposition}\label{prop:ChebyshevPullback}
Let $A\in\Z^{(n-1)\times n}$ have rank $n-1$, let $v$ be as in Theorem~\ref{thm:LissHypersurfaces} with $v_j\neq0$ for all $j$, and let $b\in\C^n$. 
Let $A_n\in\Z^{(n-1)\times n}$ be the reduced incidence matrix of the cycle graph $C_n$, whose kernel is generated by $u=(1,\dotsc,1)$.
Choose $b'\in\C^n$ with $u\cdot b'=v\cdot b$. 
Then
\[
h_{A,b}(x_1,\dotsc,x_n)\,=\,h_{A_n,b'}\big(T_{|v_1|}(x_1),\dotsc,T_{|v_n|}(x_n)\big).
\]
In particular, writing $\varphi=(T_{|v_1|},\dotsc,T_{|v_n|})\colon\C^n\to\C^n$, we have
${\cal L}_{A,b}=\varphi^{-1}({\cal L}_{A_n,b'})$ and, for $b=\mathbf 0$, $\Gamma_A=c\,\Gamma_n\circ\varphi$ for some $c\in\C^*$, where $\Gamma_n$ is the $n$-th cycle
polynomial from \cite[Def.~3.2]{Laurent2026AlgebraicBoundary}.
\end{proposition}
\begin{proof}
For $\epsilon_1\in\{\pm1\}$ we have $z_{\epsilon_1\bf 0}^v=T_{|v_1|}(x_1)\pm\sqrt{T_{|v_1|}(x_1)^2-1}$, since $z_{\epsilon_1\bf 0}^vz_{-\epsilon_1\bf 0}^v=1$ and the sum of the two is $2T_{|v_1|}(x_1)$ by Lemma~\ref{lemma:ChebyshevPoly}. 
Therefore, coordinate by coordinate and using $z_{\delta,\delta'}^v=z_{\delta,\azero}^v\, z_{\azero,\delta'}^v$, with $\beta'^{\mathbf 1}=\beta^v$, we get
\[
h_{A,b}=\prod_{\epsilon\in\{\pm1\}^n}(z_\epsilon^v-\beta^v)=\prod_{\epsilon'\in\{\pm1\}^n}\left(\prod_{j=1}^n \left(T_{|v_j|}(x_j)+\epsilon'_j\sqrt{T_{|v_j|}(x_j)^2-1}\right)-\beta'^{\bf1}\right).
\]
This is exactly $h_{A_n,b'}$ evaluated at $\varphi(x)$. 
For $b=\mathbf 0$, taking $b'=\bf0$, the statement follows from $h_A=c'\Gamma_A^2$ and $h_{A_n}=c''\Gamma_n^2$.
\end{proof}

\begin{example}
Let $w=(1,1,1,2)$, so that ${\cal L}_A=\cos(\{t: w\cdot t=0\})\subset\C^4$.
Splitting $\operatorname{supp}(w)=\{1,2\}\sqcup\{3,4\}$, we take $\ker_{\Z}(A_1)=\langle(1,1,0,0,-1)\rangle$
and $\ker_{\Z}(A_2)=\langle(0,0,1,2,1)\rangle$.
Then $\Gamma_{A_1}=x_1^2+x_2^2+x_5^2-2x_1x_2x_5-1$ and
$\Gamma_{A_2}=(x_5-x_3T_2(x_4))^2-4(1-x_3^2)(1-x_4^2)x_4^2$, and
$\Gamma_A=\operatorname{Res}_{x_5}(\Gamma_{A_1},\Gamma_{A_2})$ is irreducible of degree $9$
($4$ in $x_1,x_2,x_3$ and $8$ in $x_4$), with $37$ terms.
As predicted by Proposition~\ref{prop:ChebyshevPullback}, $\Gamma_A=\Gamma_4(x_1,x_2,x_3,T_2(x_4))$,
with $\Gamma_4$ the $4$-th cycle polynomial.
\end{example}

\section{Singular locus of Lissajous varieties}\label{sec:singularities}

In this section we describe the singular locus of Lissajous varieties by means of a sandwich-type double inclusion (Theorem~\ref{thm:SingSandwich}), expressed in terms of the component-wise cosine of subsets of $L_{A,b}$.
More precisely, we establish a chain of inclusions $\cos(Y)\subseteq{\rm Sing}(\LAb)\subseteq \cos(X)$ for suitable subsets $X,Y\subseteq L_{A,b}$.

Before turning to the singular locus, we characterize which complex affine subspaces are mapped onto the same Lissajous variety by the component-wise cosine map.
Indeed, $\LAb$ is defined through the affine subspace $L_{A,b}\subset\C^n$, but this correspondence is not one-to-one: different affine subspaces may give rise to the same Lissajous variety.
The following result makes this precise and generalizes \cite[Lemma 2.6]{Laurent2026AlgebraicBoundary}.

\begin{lemma}\label{Lemma:cos^-1}
Consider a matrix $A\in\Z^{d\times n}$ of rank $d$ and $b\in \C^n$. 
Let $\{v_j\}_{j=1}^{n-d}\subseteq\Z^n$ be a basis of the lattice $\ker(A)\cap \Z^n = \{z \in \mathbb{Z}^n \, : \, A \, z = 0 \}$.
For the Lissajous variety $\LAb = \cos (L_{A,b})$, we have $\LAb=\cos(L_{A,b}^{(s,k)})$ for any affine space 
$L_{A,b}^{(s,k)}\subseteq \C^n$ of the form  
\[
L_{A,b}^{(s,k)}=\{t\in \C^n: (s\circ v_j)\cdot t= 2\pi k_j - \tfrac{\pi}{2}b\cdot v_j \text{ for all }j=1,\dotsc,n-d\},
\]
where $s\in \{\pm 1\}^n$ and $k\in \Z^{n-d}$.
Moreover, we have the equality 
\[
\cos^{-1}(\LAb) = \bigcup_{s\in \{\pm 1\}^n,k\in \Z^{n-d}} L_{A,b}^{(s,k)}.
\]
\end{lemma}

\begin{proof}
Let $s\in\{\pm 1\}^n$ and $k\in \Z^{n-d}$, and select $x\in \cos(L_{A,b}^{(s,k)})$. 
Then,  $x=\cos(t)$  for some $t\in\C^n$ with $(s\circ v_j)\cdot t=2k_j\pi-\tfrac{\pi}{2}b\cdot v_j$ for all $j=1,\dotsc,n-d$. 
Choose $z\in\Z^n$ such that $k_j=v_j\cdot z$ for $j=1,\ldots,n-d$, and define $t'=s\circ t-2\pi z$.
With this choice we obtain $v_j\cdot t'= (v_j\circ s)\cdot t- 2\pi v_j\cdot z  = -\tfrac{\pi}{2}b\cdot v_j$, for $j=1,\dotsc,n-d$, showing $t'\in L_{A,b}$. As $\cos(t')=\cos(t)=x$, then $x\in \LAb$.
The inclusion $\LAb\subseteq \cos(L_{A,b}^{(s,k)})$ is proved analogously.

It remains to prove that  $\cos^{-1}(\LAb) \subseteq \bigcup_{s,k} L_{A,b}^{(s,k)}$.
Let $t\in \C^n$ be such that $\cos(t)\in\LAb$. 
By definition, we have $\cos(t)=\cos(t')$, for some $t' \in L_{A,b}$.
This means that there exist some $s\in \{\pm 1\}^n$ and $k'\in \Z^n$ such that $t=s\circ t'+2k'\pi$. 
Set $k=(v_1\cdot (s\circ k'),\dotsc,v_{n-d}\cdot (s\circ k'))$. 
Then, $(s\circ v_j)\cdot t=2\pi k_j-\tfrac{\pi}{2}b\cdot v_j$ for $j=1,\dotsc,n-d$, showing that $t\in L_{A,b}^{(s,k)}$.
\end{proof}

From the proof of \cite[Proposition 5.6]{BelAfiaMeroniTelen2025Chebyshev}, the singular locus of the Lissajous variety $\LAb$ is contained in the union of the rank-deficiency points for the Jacobian of the trigonometric parametrization $\phi_{A,b}$ defined as in (\ref{eq:parametrization}), and its self-intersections. 
To address the first set, define 
\begin{equation}
    S_{A,b}^{(J,h)} = \{ t\in L_{A,b}\, :\, t_j = h_j\pi, \text{ for all } j\notin J\},
\end{equation}
where $J\in {\cal I}=\{J\subset [n]\, :\, \rank(A_J)<d\}$, and $h\in \mathbb{Z}^{n-|J|}$.
A point $x\in \LAb$ makes the Jacobian of $\phi_{A,b}$ rank deficient, namely, $\rank (J_{\phi_{A,b}}(u))<d$ for some $u\in \phi_{A,b}^{-1}(x)$, if and only if $x \in \cos(S_{A,b}^{(J,h)})$. 
Indeed, the Jacobian
\[
J_{\phi_{A,b}}(u)^\top=-\left(a_1\sin(a_1\cdot u-\tfrac{\pi}{2}b_1),\dotsc,a_n\sin(a_n\cdot u-\tfrac{\pi}{2}b_n)\right)
\]
is rank deficient if and only if $\sin(t_j)=0$ for all $j\notin J$, for some $J\in {\cal I}$ and $t=A^\top u - \tfrac{\pi}{2}b\in L_{A,b}$.
Equivalently, there exist some $h\in \mathbb{Z}^{n-|J|}$ and $J\in {\cal I}$ such that $t_j=h_j\pi$.

If $x$ is a self-intersection point, then there exist $t,t'\in L_{A,b}$ with $\cos(t)=\cos(t')=x$ that give different branches at $x$. 
As $\cos(t)=\cos(t')$, there exist $s\in \{\pm 1 \}^n$ and $m\in \mathbb{Z}^n$ such that $t'=s\circ t+2\pi m$. 
As in Lemma \ref{Lemma:cos^-1}, this puts $t'\in L_{A,b}^{(s,k)}$, for a suitable $k\in\mathbb{Z}^{n-d}$, so $x\in \cos(L_{A,b}\cap L_{A,b}^{(s,k)})$. 
Moreover, for $x$ to be an actual self-intersection point, we need to have $s\circ \Row(A)\neq \Row(A)$.
To simplify the notation, let
\[
\Sigma_A=\{\sigma\in\{\pm1\}^n \, : \, \sigma\circ\Row(A)=\Row(A)\},
\]
and assume by contradiction that $s\in \Sigma_A$.
For arbitrarily small $\omega\in \Row (A)$, we can write $\cos(t+\omega)=\cos(s\circ t'+\omega)=\cos(t'+s\circ\omega)$. 
Since $s\circ\Row(A) = \Row(A)$, we have that $s\circ \omega$ is again an arbitrarily small vector of $\Row (A)$.
Hence, locally around $x$, $t$ and $t'$ parametrize the same branch. 
Hence, $s\notin \Sigma_A$.

\begin{figure}
\centering
\includegraphics[height = 6cm]{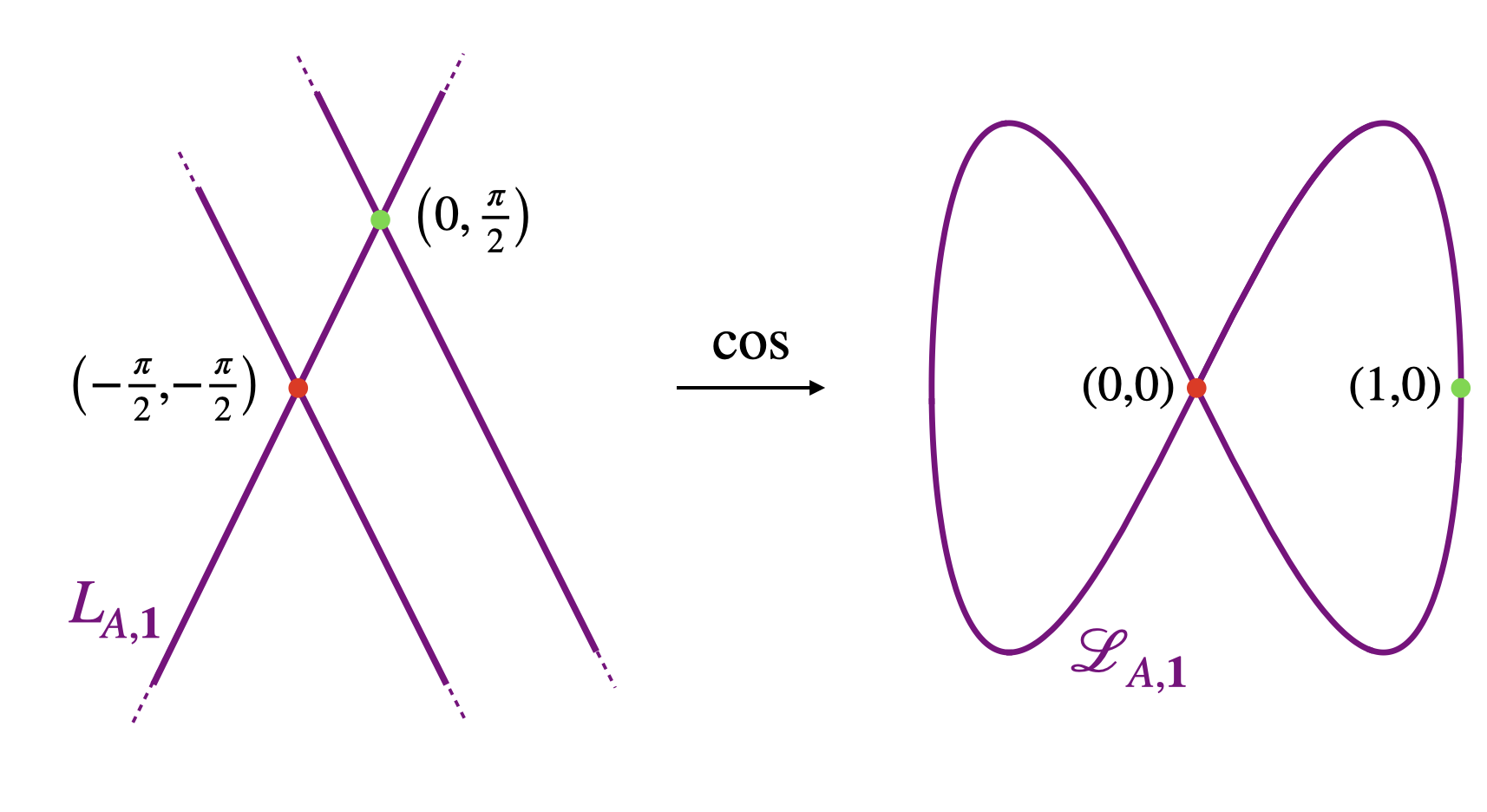}
\caption{The singular locus of the lemniscate from Example \ref{ex:ex1.1} is given by the node $\{(0,0)\}$. The line $L_{A,{\bf 1}}=\{2x-y=-\frac{\pi}{2}\}$, for $A=(1\ 2)$, intersects the lines $L^{(s,k)}_{A,{\bf 1}}=\{-2x-y=2\pi k-\frac{\pi}{2}\}$, $s=(-,+)$, for $k=0,1$, at $(0,\frac{\pi}{2})$ (green) and $(-\frac{\pi}{2},-\frac{\pi}{2})$ (red), respectively.
The intersection point $(-\frac{\pi}{2},-\frac{\pi}{2})$ is mapped to the node $(0,0)$, whereas $(0,\frac{\pi}{2})$ is mapped to the smooth point $(1,0)$: there the sign change of the first coordinate has no effect, since $t_1=0$.}
\label{fig:SingLemniscate}
\end{figure}

From the above discussion, we derive the following proposition.
\begin{proposition}\label{prop:LissSing}
With the above notation, the following inclusion holds:
\begin{equation}\label{eq:naiveSingLAb}
{\rm Sing}(\LAb)\,\subseteq\, \bigcup_{J,h} \cos(S_{A,b}^{(J,h)})\cup \bigcup_{s,k} \cos(L_{A,b}\cap L_{A,b}^{(s,k)}),
\end{equation}
where the unions range over $J\in\mathcal I$, $h\in\mathbb Z^{n-|J|}$, $s\notin\Sigma_A$, and $k\in\mathbb Z^{n-d}$.
\end{proposition}
The inclusion \eqref{eq:naiveSingLAb} is in general strict. 
For the lemniscate of Example \ref{ex:ex1.1}, with $A=\left(\begin{smallmatrix}1&2\end{smallmatrix}\right)$, $b=\bf 1$, and $s=(-,+)\notin\Sigma_A$, the set $L_{A,b}\cap L^{(s,1)}_{A,b}=\{(-\frac{\pi}{2},-\frac{\pi}{2})\}$ is mapped to the singular node $(0,0)$, while $L_{A,b}\cap L^{(s,0)}_{A,b}=\{(0,\frac{\pi}{2})\}$ is mapped to the smooth point $(1,0)$, as shown in Figure~\ref{fig:SingLemniscate}. 
In the second case the sign change acts trivially, because the flipped coordinate satisfies $\sin t_1=0$. 

Our next goal is to tighten the inclusion by only taking into account the signs of $s$ at the coordinates where the sine does not vanish.
To do so, write $W_{A,b}^{(s,k)}=L_{A,b}\cap L_{A,b}^{(s,k)}$.
For $W\subset \C^n$, define 
\[
P_W = \{i\in[n]\, :\, \exists t\in W \text{ such that } \sin(t_i)\neq 0\}.
\]
The main theorem of this section is the following.
\begin{theorem}\label{thm:SingSandwich}
    Let $A\in\Z^{d\times n}$ have rank $d$, and $b\in \C^n$. 
    With the above notation, the singular locus of the Lissajous variety $\LAb = \cos (L_{A,b})$ satisfies the following chain of inclusions: 
    \begin{equation}
        \bigcup_{s,k} \cos(W_{A,b}^{(s,k)})\,\subseteq\,{\rm Sing}(\LAb)\,\subseteq\, \bigcup_{J,h} \cos(S_{A,b}^{(J,h)})\cup \bigcup_{s,k} \cos(W_{A,b}^{(s,k)}),
    \end{equation}
    where the first union on the right-hand side ranges over $J\in\mathcal I$, $h\in\mathbb Z^{n-|J|}$, while the other two range over $W=W_{A,b}^{(s,k)}$, with $s\in \{\pm1\}^n$ and $k\in\mathbb Z^{n-d}$, satisfying the following two properties:
    \begin{itemize}
        \item[${\rm P1)}$] There is no $\sigma \in \Sigma_A$ such that $s|_{P_W}=\sigma|_{P_W}$;
        \item[${\rm P2)}$] $\rank (A_{P_W})=d$.
    \end{itemize}
\end{theorem}
Before proving Theorem \ref{thm:SingSandwich}, we need two preliminary results. 
Throughout this section, denote with $f_{A,b}=\cos|_{L_{A,b}}\colon L_{A,b}\rightarrow \LAb$. 
The first lemma states that $f_{A,b}$ is a proper map, modulo periodicity.
\begin{lemma}\label{lemma:fpropermod2pi}
    Let $(A,b)$ as before, and write $f_{A,b}=\cos|_{L_{A,b}}\colon L_{A,b}\rightarrow \LAb$.
    Suppose that the sequence $y_m=f_{A,b}(\rho_m)\subset \LAb$, with $\rho_m=A^\top \theta_m-\tfrac{\pi}{2}b\in L_{A,b}$, converges to $x\in \LAb$. 
    Then, after passing to a subsequence, there exist $\tau\in f_{A,b}^{-1}(x)$ and $(\lambda_m)_m\subset 2\pi(\Row A\cap \mathbb{Z}^n)$, such that $\rho_m-\lambda_m$ converges to $\tau$.
\end{lemma}
\begin{proof}
    Without loss of generality, we can put $A$ in the form $A=\left( D\,|\,Q\right)\in\mathbb{Z}^{d\times n}$, where $D={\rm diag}(a_{11},\dotsc,a_{dd})\in\mathbb{Z}^{d\times d}$, and $a_{jj}\neq0$. 
    
    Write $y_{mj}=(y_m)_j=\tfrac{1}{2}(\beta_jz_{mj}+\beta_j^{-1}z_{mj}^{-1})$, where $\beta_j=e^{-i\tfrac{\pi}{2}b_j}$ and $z_{mj}={\exp}({a_j\cdot i\theta_m})$.
    As $(y_m)_m$ converges, then it is bounded. 
    Hence, there exists some finite $C\geq 0$ such that $|y_{mj}|\leq C$ for all $m,j$. 
    As $\beta_j^2z_{mj}^2-2y_{mj}\beta_jz_{mj}+1=0$, we get $0<(C+\sqrt{C^2+1})^{-1}\leq|\beta_jz_{mj}|\leq C+\sqrt{C^2+1}$, meaning that also $(z_m)_m$ is bounded.
    For $j=1,\dotsc,d$, $z_{mj}=\exp(a_{jj}\cdot i \theta_{mj})$. 
    Hence, $|z_{mj}|=\exp(-a_{jj}{\rm Im}(\theta_{mj}))$, or, equivalently ${\rm Im}(\theta_{mj})=-\frac{1}{a_{jj}}\log(|z_{mj}|)$.
    From $(z_m)_m$ being bounded we can derive that also the imaginary part of the sequence $(\theta_{m})_m$ is bounded.

    Now select $k_m \in \mathbb{Z}^d$ so that the real part of $\tilde{\theta}_m=\theta_m-2\pi k_m$ lies in the box $[0,2\pi]^d$.
    As both the imaginary and the real parts of $\tilde{\theta}_m$ are bounded, then the sequence is bounded. 
    Hence, after passing to a subsequence, it converges to some $\tilde{\theta}\in\C^d$.

    Set $\lambda_m=2\pi A^\top k_m\in 2\pi(\Row(A)\cap \Z^n)$. 
    Then, $\rho_m-\lambda_m=A^\top\tilde{\theta}_m-\frac{\pi}{2}b$ converges to $A^\top\tilde{\theta}-\frac{\pi}{2}b=:\tau$. 

    Finally, $f_{A,b}(\rho_m-\lambda_m)=f_{A,b}(\rho_m)=y_m$, which converges to $x$. 
    Since $f_{A,b}(\rho_m-\lambda_m)$ also converges to $f_{A,b}(\tau)$, we obtain that $x=f_{A,b}(\tau)$.
\end{proof}
As a direct consequence of \ref{lemma:fpropermod2pi}, we get the following result.
\begin{lemma}\label{lemma:branchLAb}
    Let $x\in \LAb$ and $\epsilon>0$.
    Then there exists a neighborhood $V_\epsilon\subset \C^n$ of $x$ such that for every $y\in V_\epsilon\cap\LAb$, there exist $\rho\in f_{A,b}^{-1}(y)$ and $\tau\in f_{A,b}^{-1}(x)$ such that $|\rho-\tau|<\epsilon$.
\end{lemma}
\begin{proof}
    By contradiction, suppose that for every neighborhood $V_\epsilon$ of $x$, there exists $y\in V_\epsilon\cap\LAb$ such that for all $\rho\in f_{A,b}^{-1}(y)$ and $\tau\in f_{A,b}^{-1}(x)$, we have $|\rho-\tau|\geq\epsilon$.
    Choose $V_\epsilon=B_{\frac{1}{m}}(x)$, for $m\geq 1$.
    Then, for every $m\geq 1$, there is some $y_m\in \LAb$, with $|y_m-x|<\frac{1}{m}$, such that for every $\rho_m\in f_{A,b}^{-1}(y_m)$ and $\tau\in f_{A,b}^{-1}(x)$, we have $|\rho_m-\tau|\geq\epsilon$.
    By Lemma \ref{lemma:fpropermod2pi}, up to taking a subsequence, there exist $\tau\in f_{A,b}^{-1}(x)$ and $(\lambda_m)_m\subset 2\pi(\Row A\cap \mathbb{Z}^n)$, such that $\rho_m-\lambda_m$ converges to $\tau$.
    Set $\tau_m=\tau+\lambda_m\in L_{A,b}$.
    Then, $\tau_m\in f_{A,b}^{-1}(x)$ for every $m\geq 1$, and $|\rho_m-\tau_m|=|\rho_m-\lambda_m-\tau|\xrightarrow[m\rightarrow\infty]{} 0$.
    Hence, for $m$ large enough, $|\rho_m-\tau_m|<\epsilon$, reaching a contradiction.
\end{proof}
We are now ready to prove Theorem \ref{thm:SingSandwich}.
\begin{proof}[Proof of Theorem \ref{thm:SingSandwich}]
    We first prove that ${\rm Sing}(\LAb)\subseteq \bigcup_{J,h} \cos(S_{A,b}^{(J,h)})\cup \bigcup_{s,k} \cos(W_{A,b}^{(s,k)})$. 
    We prove it by showing that if $x$ does not belong to the right-hand side of the inclusion, then it is a regular point.
    Let $P=\{i\in[n]\, :\, x_i\neq \pm 1\}$.
    As $x\notin \bigcup_{J,h} \cos(S_{A,b}^{(J,h)})$, then $\rank(A_P)=d$ (otherwise, $x\in \cos(S_{A,b}^{(P,h)})$, for some $h\in\Z^{n-|P|}$.
    In particular, $\rank(df|_\tau(a_1),\dotsc,df|_\tau(a_d))=d$, at every $\tau\in f^{-1}(x)$, where we write $f=f_{A,b}$ to simplify the notation.
    This means that $f$ is an immersion at every $\tau\in f^{-1}(x)$.
    Select $t\in f^{-1}(x)$. 
    Since $f$ is immersive at $t$, there exists a ball $B_\epsilon(t)$, for some $\epsilon>0$, such that $f|_{L_{A,b}\cap B_\epsilon(t)}$ is an embedding. 
    In particular, ${\cal M}=f(L_{A,b}\cap B_\epsilon(t))$ is a $d$-dimensional smooth submanifold.
    By Lemma \ref{lemma:branchLAb}, applied to $x$ and $\epsilon$, there exists $V_\epsilon$ neighborhood of $x$ such that for every $y\in V_\epsilon\cap\LAb$, there exist $\rho\in f_{A,b}^{-1}(y)$ and $\tau\in f_{A,b}^{-1}(x)$ such that $|\rho-\tau|<\epsilon$.
    That is,
    \[
    V_\epsilon\cap\LAb\,\subseteq\, \bigcup_{\tau\in f^{-1}(x)}f( B_\epsilon(\tau)\cap L_{A,b}).
    \]
    Our goal is to prove that the right-hand side is equal to $V_\epsilon\cap\cal M$.
    Let $t'\in f^{-1}(x)$. 
    Since $\cos(t)=\cos(t')=x$, there exist $s\in\{\pm1\}^n$ and $m\in \Z^n$ such that $t'=s\circ t+2\pi m$.
    Let $k=-(v_1\cdot m, \dotsc, v_{n-d}\cdot m)$, with $\{v_j\}_{j=1}^{n-d}\subseteq\Z^n$ a basis of the lattice $\ker(A)\cap \Z^n$.
    From $v_j\cdot t'=-{\pi \over 2}v_j\cdot b$, we derive $(s\circ v_j)\cdot t=2\pi k_j-{\pi \over 2}v_j\cdot b$, for all $j=1,\dotsc,n-d$.
    This means that $t\in W_{A,b}^{(s,k)}$, and hence $x\in\cos(W_{A,b}^{(s,k)})$.
    For $W=W_{A,b}^{(s,k)}$ as just defined, we get $P\subset P_W$.
    Indeed, if $i \in P$, then $x_i\neq \pm 1$, and hence $\sin(t_i)\neq 0$. Since $t\in W$, then $i\in P_W$.
    In particular, since $\rank(A_P)=d$, then also $\rank A_{P_W}=d$, which settles property (i).
    Since $x\notin \cos(W)$ for $W$ satisfying both properties (P1) and (P2), then (P1) must fail.
    Namely, there exists $\sigma \in \Sigma_A$ such that $s|_{P_W}=\sigma|_{P_W}$.
    In particular, $s|_P=\sigma|_P$.
    Since $t'=s\circ t+2\pi m$, one can verify that $t'=\sigma\circ t+2\pi \lambda$, for some $\lambda\in \Z^n$.
    Let $\omega\in \Row(A)$.
    Then, $f(t+\omega)=f(\sigma\circ t'+\omega)=f(t'+\sigma\circ\omega)$. 
    As $\sigma \in \Sigma_A$, we have that $\omega'=\sigma\circ\omega\in \Row (A)$.
    Moreover, $||\omega||<\epsilon$ if and only if $||\omega'||<\epsilon$, because $\omega\mapsto\sigma\circ\omega$ is an isometry.
    This shows that $f( B_\epsilon(t')\cap L_{A,b})=f( B_\epsilon(t)\cap L_{A,b})=\cal M$, and, iterating for all $\tau \in f^{-1}(x)$, that $\bigcup_{\tau\in f^{-1}(x)}f( B_\epsilon(\tau)\cap L_{A,b})=\cal M$.
    Then, $V_\epsilon\cap\LAb\subseteq{\cal M}\subseteq\LAb$. 
    This means that $V_\epsilon\cap\LAb=V_\epsilon\cap{\cal M}\ni x$, and hence that $x$ is regular.

    We now prove the second inclusion, showing that $\cos(W_{A,b}^{(s,k)})\subseteq{\rm Sing}(\LAb)$, for all $W_{A,b}^{(s,k)}$ satisfying (P1) and (P2).
    Let $t'\in W_{A,b}^{(s,k)}=W$ such that $\sin(t_i')\neq 0$ for all $i\in P_W$. 
    Denote $t=s\circ t'-2\pi h$, where $h$ satisfies $v_j\cdot h= k_j$ for all $j=1,\dotsc,n-d$.
    It is clear that $\cos(t)=\cos(t')=x$ and one can check that $t\in L_{A,b}$.
    We claim that $t$ and $t'$ represent different branches through $x$.
    Notice that $t$ and $t'$ share the same set of indices where the sine is non-zero, namely $P_W$, as $\sin(t)=s\circ\sin(t')$.
    Then, by property (P2), $f=f_{A,b}$ is an immersion at both $t$ and $t'$, and hence locally an embedding.
    This means that there exists some $\epsilon>0$ such that $f(L_{A,b}\cap B_\epsilon(t))=\cal M$ and $f(L_{A,b}\cap B_\epsilon(t'))=\cal M'$ are smooth $d$-dimensional submanifold.

    Suppose, by contradiction, that $x$ is a regular point of $\LAb$.
    Since $\LAb$ is irreducible of dimension $d$, there is a neighbourhood $U$ of $x$ such that $\LAb\cap U$ is a $d$-dimensional manifold. 
    For $\epsilon$ small, $\mathcal M,\mathcal M'\subset U$ are $d$-dimensional submanifolds of $\LAb\cap U$, hence open in it, and both contain $x$.
    Thus for $\epsilon_1\ll\epsilon$ we have $\mathcal M_{\epsilon_1}\subseteq\mathcal M'_{\epsilon}$, where $\mathcal M_{\epsilon_1}=f(L_{A,b}\cap B_{\epsilon_1}(t))$.
    That is, for every $\omega\in\Row(A)$ with $\|\omega\|<\epsilon_1$ there is $\omega'\in\Row(A)$ with $\|\omega'\|<\epsilon$ such that $f(t+\omega)=f(t'+\omega')$.

    Reasoning component-wise and using the fact that $\omega,\omega'\approx 0$ (for $\epsilon$ very small), one shows~that
    \[
    \omega'=\sigma \circ \omega, 
    \text{ where }
    \begin{cases}
        \sigma_i=s_i \text{ for all } i\in P_W,\\
        \sigma_i \in \{\pm 1\} \text{ otherwise}.
    \end{cases}
    \]
    That is, there exists $\sigma \in \{\pm 1\}^n$ such that $\omega'=\sigma\circ \omega$ and $\sigma|_{P_W}=s|_{P_W}$.
    Set $E_\sigma = \{\omega\in \Row(A)\, :\, \sigma\circ \omega\in \Row(A)\}$, which is a linear subspace of $\Row(A)$.
    We have just shown that for all $\omega \in \Row(A)$ of size $||\omega||<\epsilon_1$, there exists $\sigma\in \{\pm 1\}^n$ with $\sigma|_{P_W}=s|_{P_W}$ such that $\omega\in E_\sigma$.
    Notice that there are only finitely many $E_\sigma$'s (namely, $2^{n-|P_W|}$) and $\Row(A)\cap B_{\epsilon_1}(0)\subset\bigcup_\sigma E_\sigma$.
    This means that there exists $\tilde{\sigma}$ such that $E_{\tilde{\sigma}}=\Row(A)$.
    In particular, we have found $\tilde{\sigma}\in \{\pm 1\}^n$ such that $\tilde{\sigma}\circ \Row(A)=\Row(A)$ and $\tilde{\sigma}|_{P_W}=s|_{P_W}$.
    As $\tilde{\sigma}\in \Sigma_A$, this contradicts (P1). 
    Hence, $x\in {\rm Sing}(\LAb)$.

    This shows $\cos(W^\circ)\subseteq{\rm Sing}(\LAb)$ with $W^\circ:=\{t'\in W:\sin t'_i\neq0\ \forall i\in P_W\}$.
    The set $W^\circ$ is dense in $W$, the map $\cos$ is continuous and ${\rm Sing}(\LAb)$ is closed, hence $\cos(W)\subseteq\overline{\cos(W^\circ)}\subseteq{\rm Sing}(\LAb)$.
\end{proof}

\begin{remark}
    The singular locus ${\rm Sing}(\LAb)$ is contained in a union of Lissajous varieties.
    Indeed, each $S_{A,b}^{(J,h)}$ and each $W_{A,b}^{(s,k)}$ is an affine subspace of $\C^n$ whose direction space is defined over $\mathbb Q$ (it is $\Row(A)$ intersected with $\{\tau_j=0,\ j\notin J\}$, respectively with $(s\circ\ker A)^{\perp}$).
    Such a subspace has the form $\Row(A')-\frac{\pi}{2}b'$ for some integer matrix $A'$ of full row rank (possibly empty, when the subspace is a point) and some $b'\in\C^n$, so its image under $\cos$ is the Lissajous variety $\mathcal L_{A',b'}$.
    Notice that there are only finitely many distinct images of $S_{A,b}^{(J,h)}$ and $W_{A,b}^{(s,k)}$ under the map $\cos$, so only finitely many such varieties occur.
    The same holds for the union on the left-hand side of Theorem \ref{thm:SingSandwich}, which is contained in ${\rm Sing}(\LAb)$.
\end{remark}

\subsection*{Singular locus of Lissajous curves}
In the case when $d=1$, namely of Lissajous curves, the picture becomes particularly easy to describe.
Indeed, as we will see in a moment, the sets $S_{A,b}^{(J,h)}$ play no role in the description of the singularities, and the sets $W_{A,b}^{(s,k)}$ are single points (or empty).

\begin{corollary}\label{cor:Singcurves}
    Let $A\in\Z^{1\times n}$ be non-zero, and $b\in \C^n$. 
    With the notation of Theorem \ref{thm:SingSandwich}, the singular locus of the Lissajous curve $\LAb = \cos (L_{A,b})$ is characterized by: 
    \begin{equation}
        {\rm Sing}(\LAb) = \bigcup_{s,k} \cos(W_{A,b}^{(s,k)}),
    \end{equation}
    where the union ranges over the $(s,k)$ for which $W_{A,b}^{(s,k)}$ satisfies (P1).
\end{corollary}

\begin{proof}
    First, we show that every point in $\cos(S_{A,b}^{(J,h)})$ is smooth. 
    Assume for the moment that $A=(a_1,\dotsc,a_n)$ has no zero entry.
    This mean that ${\cal I}=\{\emptyset\}$ and hence that $S_{A,b}^{J,h}=S_{A,b}^{(\emptyset,h)}$.
    Let $x\in\cos( S_{A,b}^{(\emptyset,h)})$. Then, $x_j=(-1)^{h_j}$, for all $j=1,\dotsc,n$.
    Select $t\in f^{-1}(x)$, where $f=f_{A,b}$ as usual.
    As $x_j=(-1)^{h_j}$, then $t_j=\tilde{h}_j\pi$ with $(-1)^{\tilde{h}_j}=(-1)^{h_j}$.
    
    We now perturb $t\in L_{A,b}$ along $L_{A,b}$.
    For any $\epsilon>0$ small, consider $\omega\in\C$, with $|\omega|<\epsilon$.
    Write $a=(a_1,\dotsc,a_n)^{\top}$.
    We obtain that
    \[
    \cos(t_j+\omega a_j)=\cos(\tilde{h}_j\pi+\omega a_j)=(-1)^{\tilde{h}_j}cos(\omega a_j)=x_j\cos(\omega a_j).
    \]
    Hence, for any $t\in f^{-1}(x)$, $f(B_{\epsilon}(t)\cap L_{A,b})=\{x\circ \cos(\omega a)\,:\, |\omega|<\frac{\epsilon}{||a||}\}=:\cal M$.
    As in the proof of \ref{thm:SingSandwich}, by Lemma \ref{lemma:branchLAb}, for every $\epsilon$, there exists a neighborhood $V_\epsilon$ of $x$ such that $V_\epsilon\cap\LAb\subseteq \bigcup_{\tau\in f^{-1}(x)}f( B_\epsilon(\tau)\cap L_{A,b})=\cal M$.

    We now prove that ${\cal M} \cap V_\epsilon$ for some $\epsilon>0$ small is smooth. 
    Indeed, if we expand in Taylor series around $\omega=0$, we get 
    \[
    x_j\cos(\omega a_j)=x_j-\frac{x_j a_j^2}{2}\omega^2+O(\omega^4).
    \]
    Substituting $\omega^2=z$ and using the fact that the function is even, we can re-parametrize ${\cal M}=G(B_{(\epsilon/||a||)^2}(0))$, where $G_j(z)=x_j-\frac{x_j a_j^2}{2}z+O(z^2)$.
    In particular, $\frac{dG_j}{dz}|_{z=0}=-\frac{x_j a_j^2}{2}\neq 0$.
    Hence $G$ is an immersion at $z=0$ and locally an embedding on some disk $D_r(0)$. 
    Choosing $\epsilon>0$ so that $\frac{\epsilon^2}{||a||^2}<r$, we get that $\cal M$ is smooth.
    Since $V_\epsilon\cap\LAb=V_\epsilon\cap {\cal M}$, locally around $x$, $\LAb$ is a smooth curve.
    Hence, $x\notin{\rm Sing}(\LAb)$.
    If some entries of $A$ are zero, the proof is analogue, as $G$ remains an immersion at $z=0$.

    We now look at $\cos(W_{A,b}^{(s,k)})$, for some $(s,k)$, satisfying properties (P1) and (P2) of Theorem \ref{thm:SingSandwich}.
    We show that, for $d=1$, (P1) implies (P2).
    Indeed, $\Sigma_A=\{\sigma\in\{\pm 1\}^n\, :\, \sigma\circ a = \pm a\}$, namely it is given by all $\sigma \in \{\pm 1\}^n$ such that every $\sigma_i=+1$ (respectively, $-1$) for all $i\in [n]$ such that $a_i\neq 0$. 
    Suppose by contradiction that (P2) does not hold for some $W=W_{A,b}^{(s,k)}$. 
    Then, $\rank(A_{P_W})<1$ and hence $A_{P_W}=({\bf 0})$.
    In particular, $P_{W}=\{i\in[n]\,:\,\sin(t_i)\neq 0 \text{ for some } t\in W\}\subseteq \{i\in[n]\,:\,a_i=0\}$.
    Construct $\sigma$ in this way: if $i\in P_W$, then $\sigma_i=s_i$ and if $i\notin P_W$, then $\sigma_i=+1$.
    By construction, $\sigma\in \Sigma_A$ and $\sigma|_{P_W}=s|_{P_W}$, contradicting (P1).

    We now conclude applying Theorem \ref{thm:SingSandwich}, and obtain that 
    ${\rm Sing}(\LAb) = \bigcup_{s,k} \cos(W_{A,b}^{(s,k)})$, where each $W_{A,b}^{(s,k)}$ satisfies (P1).
\end{proof}

\begin{example}\label{ex:247SingCurves}
    Let $A=\left(\begin{smallmatrix} 2 & 4 & 7 \end{smallmatrix}\right)$, so that $L_{A,b}=\{t\,a-\tfrac{\pi}{2}b : t\in\C\}$ with $a=(2,4,7)$.
That is, $\mathcal L_{A,b}$ is the curve
$(\cos(2t-\tfrac{\pi}{2}b_1),\cos(4t-\tfrac{\pi}{2}b_2),\cos(7t-\tfrac{\pi}{2}b_3))$, $t\in \C$.
We have $\Sigma_A=\{\pm(1,1,1)\}$, and $\ker(A)\cap\mathbb Z^3$ has basis $v_1=(2,-1,0)$, $v_2=(7,0,-2)$.
The point $ta-\tfrac{\pi}{2}b$ lies in $W^{(s,k)}$ if and only if
\begin{equation}\label{eq:W247}
( s\circ v_j\cdot a)\,t \;=\; 2\pi k_j-\pi\!\!\sum_{i:\,s_i=-1}\!\!(v_j)_i\,b_i,
\quad s\circ v_1\cdot a = 4(s_1-s_2),\quad s\circ v_2\cdot a=14(s_1-s_3).
\end{equation}
By Corollary~\ref{cor:Singcurves}, $\mathrm{Sing}(\mathcal L_{A,b})$ is the union of the points $\cos(W^{(s,k)})$ over the $W^{(s,k)}$ that satisfy (P1), i.e., those for which $s$ is not constant on $P_W$.

For the case $b=\mathbf 0$, $\LA$ is shown in Figure~\ref{fig:Curve247}, on the left; here $W^{(s,k)}=W^{(-s,-k)}$, so we may assume $s_1=1$.
The signs $s=(+,-,+)$ and $s=(+,-,-)$ do not give any $W$ satisfying (P1): \eqref{eq:W247} forces $\sin 4t=0$, resp.\ $\sin 2t=0$, so that $s$ is constant on $P_W$.
For $s=(+,+,-)$, \eqref{eq:W247} gives $k_1=0$ and $W^{(s,(0,k_2))}=\{\tfrac{k_2\pi}{14}a\}$.
We have $3\in P_W$ if and only if $\sin\tfrac{k_2 \pi}{2}\neq0$, i.e., $k_2$ is odd.
Moreover $P_W\cap\{1,2\}=\emptyset$ if and only if $\tfrac{k_2\pi}{14}\in\tfrac{\pi}{2}\mathbb Z$.
Hence $W$ satisfies (P1) if and only if $k_2$ is odd and $k_2\not\equiv 7 \bmod 14$.
The resulting point only depends on $k_2$ modulo $14$ and up to sign.
Therefore the singular locus consists of three points:
\[
\mathrm{Sing}(\LA)=\big\{(\cos\tfrac{\pi k}{7},\cos\tfrac{2\pi k}{7},0) : k=1,3,5\big\}.
\]

For the case $b=\mathbf 1$, $\mathcal L_{A,\mathbf 1}$ is the curve $(\sin 2t,\sin 4t,\sin 7t)$ (see Figure~\ref{fig:Curve247}, on the right), and $s$ and $-s$ have to be treated separately.
For $s=(-,-,+)$, $(-,+,-)$, and $(+,-,-)$, \eqref{eq:W247} has no solution, so $W=\emptyset$; for $s=(+,-,+)$ and $(-,+,+)$, as in the case $b=\mathbf 0$, \eqref{eq:W247} forces $\sin(4t-\tfrac{\pi}{2})=0$, resp.\ $\sin(2t-\tfrac{\pi}{2})=0$, so $W$ does not satisfy (P1).
The sign $s=(+,+,-)$ gives $k_1=0$ and $W^{(s,(0,k_2))}=\{\tfrac{(k_2+1)\pi}{14}a-\tfrac{\pi}{2}\mathbf 1\}$.
Here $7t-\tfrac{\pi}{2}=\tfrac{k_2\pi}{2}$, so $3\in P_W$ if and only if $k_2$ is odd. For $k_2=2m-1$, i.e., $t=\tfrac{\pi m}{7}$, we have $\sin(2t-\tfrac{\pi}{2})=-\cos\tfrac{2\pi m}{7}\neq0$, so $W$ satisfies (P1) if and only if $k_2$ is odd.
The resulting point only depends on $m$ modulo $7$. 
Therefore there are seven singular points:
\[
\mathrm{Sing}(\mathcal L_{A,\mathbf 1})=\big\{(\sin\tfrac{2\pi m}{7},\sin\tfrac{4\pi m}{7},0) : m=0,\dots,6\big\}.
\]
    
    \begin{figure}
    \centering
    \includegraphics[height = 6cm]{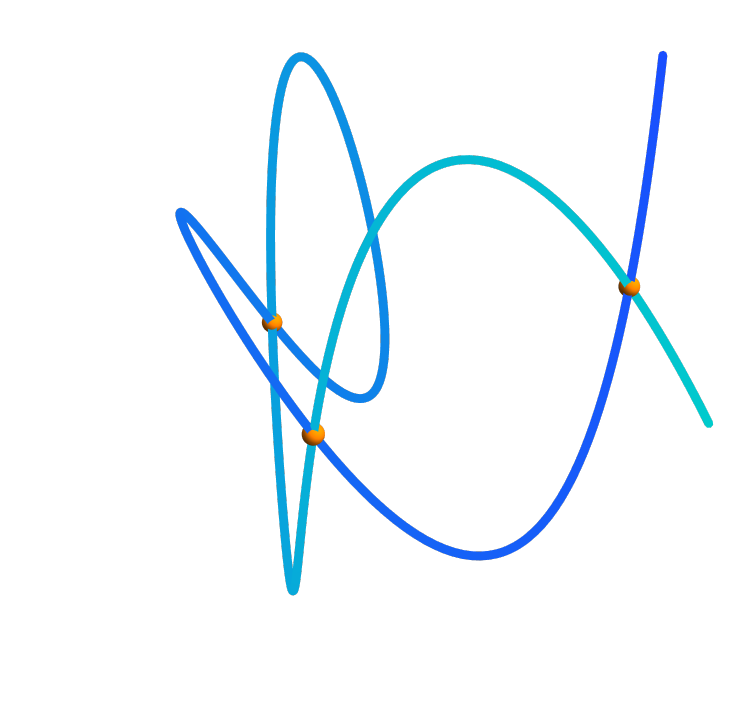} 
    \quad \quad 
    \includegraphics[height = 6cm]{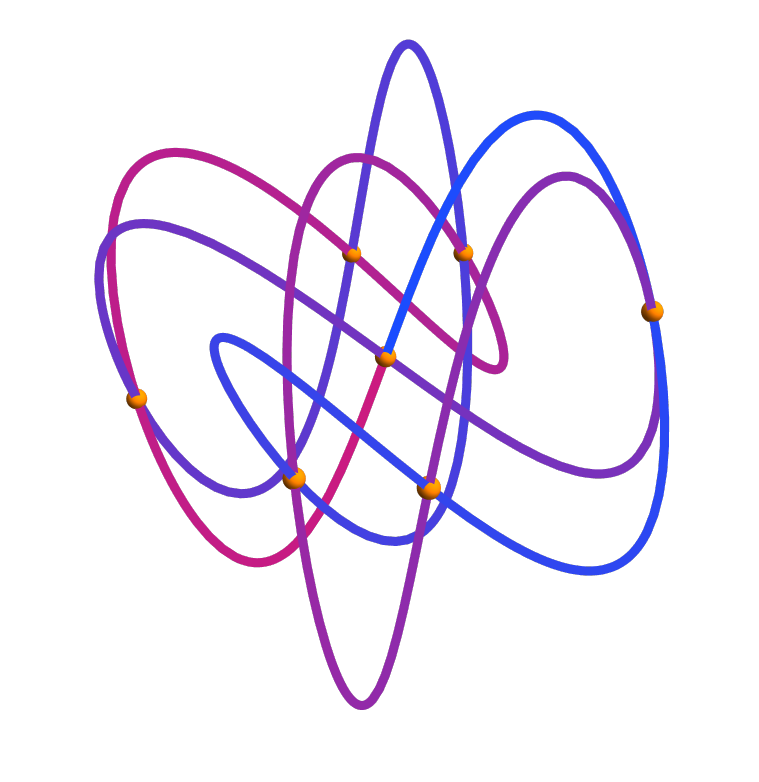}
    \caption{The Lissajous curves $\LA$ (on the left) and ${\cal L}_{A,\bf 1}$ (on the right) with $A=\left(\begin{smallmatrix} 2 & 4 & 7\end{smallmatrix}\right)$, and their singular loci (the golden bullet), from Example \ref{ex:247SingCurves}.}
    \label{fig:Curve247}
    \end{figure}
\end{example}

\begin{remark}
    We compare Corollary \ref{cor:Singcurves} with \cite[Theorem 3.3]{BelAfiaMeroniTelen2025Chebyshev}.
    Let $A=(a_1\ a_2)$ with $a_1,a_2>0$, $\gcd(a_1,a_2)=1$, and $b=\bf 0$.
    Then $\Sigma_A=\{\pm(1,1)\}$, so (P1) forces $s=\pm(1,-1)$; take $s=(1,-1)$, the other case being symmetric.
    A point $t=(a_1u,a_2u)\in L_{A,0}=\Row(A)$ lies in $W_{A,0}^{(s,k)}$ for some $k$ if and only if there exist $t'=(a_1u',a_2u')\in L_{A,0}$ and $m\in\Z^2$ with $t'=s\circ t+2\pi m$, i.e.
    \[
        a_1(u'-u)=2\pi m_1,\qquad a_2(u'+u)=2\pi m_2 .
    \]
    Eliminating $u'$ gives $u=\pi(a_1m_2-a_2m_1)/(a_1a_2)$, and since $\gcd(a_1,a_2)=1$ the numerator takes every integer value $r$.
    Hence $W_{A,0}^{(s,k)}$ is the point $t=\left({\pi r \over a_2},{\pi r\over a_1}\right)$.
    Moreover, (P1) holds if and only if $P_W=\{1,2\}$, i.e., $a_1\nmid r$ and $a_2\nmid r$.
    Therefore
    \[
        {\rm Sing}(\LA)=\Big\{\Big(\cos\tfrac{\pi r}{a_2},\cos\tfrac{\pi r}{a_1}\Big)\,:\,r\in\Z,\ a_1\nmid r,\ a_2\nmid r\Big\}.
    \]

    Let $k\in\{1,\dots,a_2-1\}$ and $l\in\{1,\dots,a_1-1\}$ be such that $k\equiv\pm r\pmod{2a_2}$ and $l\equiv\pm r\pmod{2a_1}$ (these exist since $a_2\nmid r$ and $a_1\nmid r$), so that $x=(\cos\frac{k\pi}{a_2},\cos\frac{l\pi}{a_1})$.
    Then $k\equiv l\equiv r\pmod 2$, and conversely every pair with $k\equiv l \pmod 2$ arises by the Chinese remainder theorem, since $\gcd(2a_1,2a_2)=2$.
    This is exactly the set described in \cite[Theorem 3.3]{BelAfiaMeroniTelen2025Chebyshev}.
\end{remark}

\subsection*{Singular locus of Lissajous hypersurfaces}
If $A$ is coloop-free and $n\geq 4$, then, as for curves, the sets $S_{A,b}^{(J,h)}$ are not needed to describe the singular locus of the Lissajous hypersurface $\LAb$.
\begin{corollary}\label{cor:hypersurfaces}
    Let $A\in\Z^{(n-1)\times n}$ have rank $n-1$ and ${\rm CL}_A=0$, let $n\ge4$, and let $b\in\C^n$.
    Then
    \begin{equation}
        {\rm Sing}(\LAb)=\bigcup_{(s,k)}\cos\big(W_{A,b}^{(s,k)}\big),
    \end{equation}
    where the union ranges over all $(s,k)\in\{\pm1\}^n\times\Z$ such that $W_{A,b}^{(s,k)}\neq\emptyset$,
    $|P_W|\ge n-1$ and $s|_{P_W}$ is not constant.
    Equivalently, it ranges over the $(s,k)$ satisfying (P1) and (P2) of Theorem \ref{thm:SingSandwich}.
\end{corollary}

\begin{proof}
    Let $v$ be the primitive generator of $\ker A\cap\Z^n$.
    Since ${\rm CL}_A=0$, all entries of $v$ are non-zero, so any $n-1$ columns of $A$ are linearly independent
    and $\Sigma_A=\{\pm(1,\dots,1)\}$.
    This gives the equivalence with (P1) and (P2) stated in the corollary.

    By Theorem \ref{thm:SingSandwich} it suffices to show that
    $\cos(S_{A,b}^{(J,h)})\subseteq\bigcup_{W}\cos(W)$, with the union over the $W$ that satisfies (P1) and (P2).
    Since any $n-1$ columns are independent, $J\in{\cal I}$ if and only if $|J|\le n-2$.
    Hence, if $t\in S_{A,b}^{(J,h)}$, there are two indices $i\neq j$ with $t_i=h_i\pi$ and $t_j=h_j\pi$.
    Let $s$ be the signing with $s_i=s_j=-1$ and $s_l=+1$ otherwise, and let $k=-(v_ih_i+v_jh_j)\in\Z$.
    Using $v\cdot t=-\frac{\pi}{2}v\cdot b$, we get
    $(s\circ v)\cdot t=v\cdot t-2(v_it_i+v_jt_j)=-\frac{\pi}{2}v\cdot b+2\pi k$, so $t\in W:=W_{A,b}^{(s,k)}$.

    The set $W$ is given by $v\cdot t'=-\frac{\pi}{2}v\cdot b$ and $v_it'_i+v_jt'_j=-\pi k$.
    As $n\ge4$, the remaining at least two coordinates $t'_l$ satisfy one relation with non-zero coefficients,
    so every coordinate $t'_m$ is a non-constant function on $W$. Thus $P_W=[n]$.
    Since $s$ is not constant, (P1) holds, and (P2) holds because $\rank(A)=n-1$.
    Therefore $\cos(t)\in\cos(W)$, with $W$ satisfying (P1) and (P2).
\end{proof}

\begin{remark}
    If ${\rm CL}_A>0$, the coordinates indexed by the coloops are free (since $e_j\in\Row(A)$ for each coloop $j$), so
    $\LAb=\mathcal L_{A',b'}\times\C^{{\rm CL}_A}$, where $A'$ is obtained by removing the coloops.
    Hence ${\rm Sing}(\LAb)={\rm Sing}(\mathcal L_{A',b'})\times\C^{{\rm CL}_A}$, and the corollary applies to $A'$
    whenever $n-{\rm CL}_A\ge4$.
\end{remark}

\subsection*{Acknowledgements}
I am deeply grateful to my advisor Simon Telen for his many suggestions and for checking the results, and to my advisor Monique Laurent for helpful discussions.

I used Anthropic’s Claude (Sonnet 5.5 and Opus 5.5) in the final stage of writing to identify possible gaps in the mathematical arguments, help refine some proofs, and improve the English exposition. 
I take full responsibility for the content of this paper.

This work has been supported by the European Union's HORIZON-MSCA-2023-DN-JD programme under the Horizon Europe (HORIZON) Marie Sklodowska-Curie Actions, grant agreement
101120296 (TENORS).

\bibliographystyle{unsrt}

\bibliography{bibliography.bib}

\end{document}